%% file: sv+_mc.tex
\documentclass{m2an}

\usepackage[T1]{fontenc}
\usepackage[english]{babel}
\usepackage{amssymb}
\usepackage{amsmath}
\usepackage[dvips]{epsfig}
\usepackage{color}
\usepackage{setspace}
\usepackage{graphics,graphicx,color}
\usepackage{epsfig}
\usepackage[colorlinks=true]{hyperref} %% Warning: when you first run your tex file, some errors might occur, please just
\allowdisplaybreaks

\newcommand\R{\mathbb{R}}
\graphicspath{{Figures/}}
\DeclareGraphicsExtensions{.jpg}

\def \trait (#1) (#2) (#3){\vrule width #1pt height #2pt depth #3pt}
\def \Box { \hfill
	            \trait (0.1) (5) (0)
                 	\trait (5) (0.1) (0)
            	    \kern-5pt
	            \trait (5) (5) (-4.9)
               	\trait (0.1) (5) (0)}

\begin{document}

\title[New Multilayer Saint-Venant System]{A multilayer Saint-Venant system with mass exchanges for Shallow Water flows.\\
Derivation and numerical validation}
\runningtitle{An Exchanging
Mass Multilayer Saint-Venant System }

\author{E. Audusse}
\address{Univ. Paris 13, Institut Galil\'ee,  99 avenue
Jean-Baptiste Cl\'ement 93430 Villetaneuse, France;
\email{audusse@math.univ-paris13.fr}}

\author{M.O. Bristeau}
\address{INRIA Paris-Rocquencourt, B.P.~105, 78153 Le Chesnay Cedex, France;
\email{Marie-Odile.Bristeau@inria.fr}}

\author{B. Perthame}
\sameaddress{2}\secondaddress{Lab. J.-L. Lions, Univ. P. et M. Curie, BC187, 4 place Jussieu, F75252 Paris cedex 05;
\email{benoit.perthame@upmc.fr}}

\author{J. Sainte-Marie}
\sameaddress{2}\secondaddress{Saint-Venant Laboratory, 6 quai Watier, 78400 Chatou, France;
\email{Jacques.Sainte-Marie@inria.fr}}

\thanks{This work has been achieved while the authors were involved in the ANR project METHODE (http://methode.netcipia.net)}



\begin{abstract}
The standard multilayer Saint-Venant system consists in introducing fluid
layers that are advected by the interfacial velocities. As a consequence there is no mass
exchanges between these layers and each layer is described by its height and its average
velocity.

Here we introduce another multilayer system with mass exchanges between the neighborhing
layers where the unknowns are a total height of water and an average velocity per layer.
We derive it from Navier-Stokes system with an hydrostatic pressure and prove energy and
hyperbolicity properties of the model. We also give  a kinetic interpretation leading to
effective numerical schemes with positivity and energy properties. Numerical tests show
the versatility of the approach and its ability to compute recirculation cases with wind
forcing.

%In \cite{audusse}, Audusse has proposed a
%multilayer version of the Saint-Venant system 
%derived in Perthame and Gerbeau \cite{gerbeau}. 
%The system proposed in \cite{audusse} corresponds 
%to a superposition of non miscible fluids 
%i.e. without mass exchange between the layers. 
%This limitation vanishes in the model proposed here 
%since it allows the fluid to circulate from one layer to the connected ones. 
%The main properties of the model are exhibited.
%A kinetic interpretation and some numerical simulations 
%including a recirculation case with wind forcing 
%are also given.
\end{abstract}

\maketitle

%%%%%%%%%%%%%%%%%%%%%%%%%%%%%%%%%%%%%%%%%%%%%%%%%%%%%%%%%%%%%%%%%%%%%%%%%%%%%%%%%%%%%%%%%%%%%%%
\section{Introduction}
\label{sec:intro}
%%%%%%%%%%%%%%%%%%%%%%%%%%%%%%%%%%%%%%%%%%%%%%%%%%%%%%%%%%%%%%%%%%%%%%%%%%%%%%%%%%%%%%%%%%%%%%%
Due to computational issues associated with the free surface Navier-Stokes or
Euler equations, the simulations of geophysical flows are often carried
out with shallow water type
models of reduced complexity. Indeed, for vertically
averaged models such as the Saint-Venant system \cite{saint-venant},
efficient and robust numerical techniques
(relaxation schemes \cite{bouchut}, kinetic schemes
\cite{perthame},\ldots) are available and avoid to deal with moving meshes.

\begin{figure}
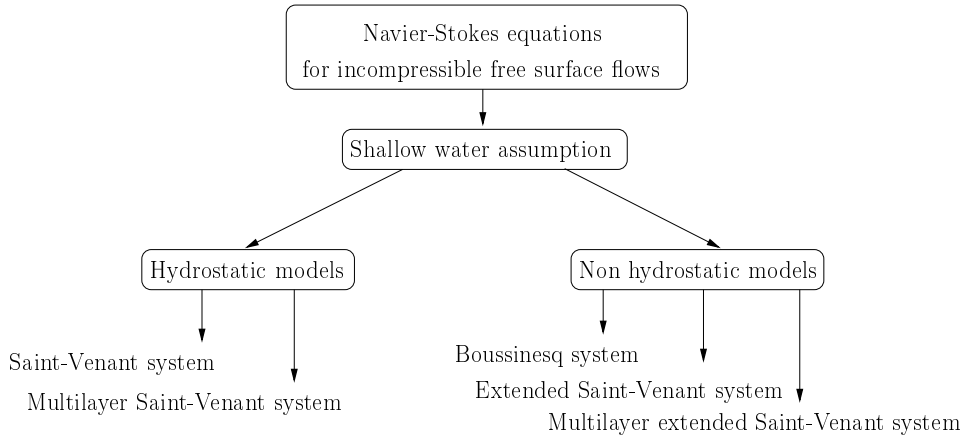

\begin{center}
\resizebox{13cm}{!}{\input Figures/models_mc.pstex_t}
%\resizebox{13cm}{!}{\input Figures/models_mc.pdf_t}
\caption{Averaged models derived from Navier-Stokes equations.}
\label{fig:shallow}
\end{center}
\end{figure}

Non-linear shallow water equations model the dynamics of a shallow,
rotating layer of homogeneous incompressible fluid and are typically
used to describe vertically averaged flows in two or three
dimensional domains, in terms of horizontal velocity and depth
variation, see Fig.~\ref{fig:shallow}.

The classical Saint-Venant system \cite{saint-venant} with viscosity
and friction \cite{decoene,saleri,gerbeau,marche} is well suited for the modeling
of dam breaks or hydraulic jumps. The extended version of the
Saint-Venant system proposed by Bristeau and Sainte-Marie \cite{JSM_DCDS} dropping the hydrostatic assumption is well adapted 
for the modeling of gravity waves propagation.

Considering flows with large friction coefficients, 
with significant water depth or with important wind effects, 
the horizontal velocity can hardly be approximated 
--~as in the Saint-Venant system~-- 
by a vertically constant velocity \cite{salencon}. To drop this limitation a multilayer Saint-Venant model is often used where each layer is described
by its own height, its own velocity and is advected by the flow (see \cite{audusse,bristeau3,bristeau2} and the
references therein). This advection property induces that there is no mass exchanges
between neighborhing layers and makes a close relation to models for two non-miscible
fluids (see \cite{bouchut1,pares,pares1}) for instance). In \cite{audusse} the multilayer strategy was formally
derived from the Navier-Stokes system with hydrostatic hypothesis departing from an
earlier work \cite{bristeau2} introducing a vertical partition of water height.

Here, we derive another and simpler multilayer model where we prescribe the vertical
discretization of the layers taking in to account the (unknown) total height of water.
Using a Galerkin approximation in lagrangian formulation, we obtain a system where the
only additional unknowns are the layers velocities.
%To drop this limitation a multilayer Saint-Venant model 
%has been introduced in \cite{audusse,bristeau3,bristeau2},
%where the different layers have their own velocities.
%This first approach was closely related to some works \cite{bouchut1,pares,pares1}
%devoted to the solution of a system of two non-miscible fluids.
%In this paper we propose a new class of multilayer Saint-Venant systems
%that allows the fluid to circulate from one layer to the connected ones. 
%In~\cite{audusse} the multilayer strategy was introduced 
%through a physical partition (along the vertical axis) of the water height.
%This leads to one continuity equation for each layer 
%and the interfaces are advected by the flow, this approach  
%corresponds to a superposition of non miscible fluids. 
%In our model the multilayer approach consists in prescribing 
%the vertical discretization of the horizontal velocity
%by performing a Galerkin approximation in Lagrangian formulation.
This leads to a global continuity equation 
and allows mass exchanges between layers.

The objective of the paper is to present the derivation 
of this new multilayer model and to
exhibit its main properties (hyperbolicity, energy equality,
\ldots). Some simulations performed with a kinetic scheme
\cite{bristeau} are presented at the end of the paper.

The paper is organized as follows. In Section \ref{sec:simplified},
we first present, in a simplified case, 
the formulation of the new multilayer Saint-Venant system
starting from the hydrostatic Euler equations. In Section \ref{sec:multi}, we recall the
Navier-Stokes system with a free moving boundary and its closure, and 
the Shallow Water system. We also introduce
the multilayer formulation in the context of the hydrostatic assumption.
In Section \ref{sec:properties} we
examine the main properties of the multilayer system and
present a kinetic interpretation of the proposed model. This kinetic formulation leads
to a numerical scheme detailed in Section \ref{sec:numeric} where some numerical
simulations are also shown.
%%%%%%%%%%%%%%%%%%%%%%%%%%%%%%%%%%%%%%%%%%%%%%%%%%%%%%%%%%%%%%%%%%%%%%%%%%%%%%%%%%%%%%%%%%%%%%%
\section{A simplified case}
\label{sec:simplified}
%%%%%%%%%%%%%%%%%%%%%%%%%%%%%%%%%%%%%%%%%%%%%%%%%%%%%%%%%%%%%%%%%%%%%%%%%%%%%%%%%%%%%%%%%%%%%%%
Before deriving the complete version of the multilayer system, 
we illustrate the approach in a simple situation. 
Moreover this case emphasizes the main differences with
the multilayer system proposed by Audusse \cite{audusse}.

\begin{figure}[h]
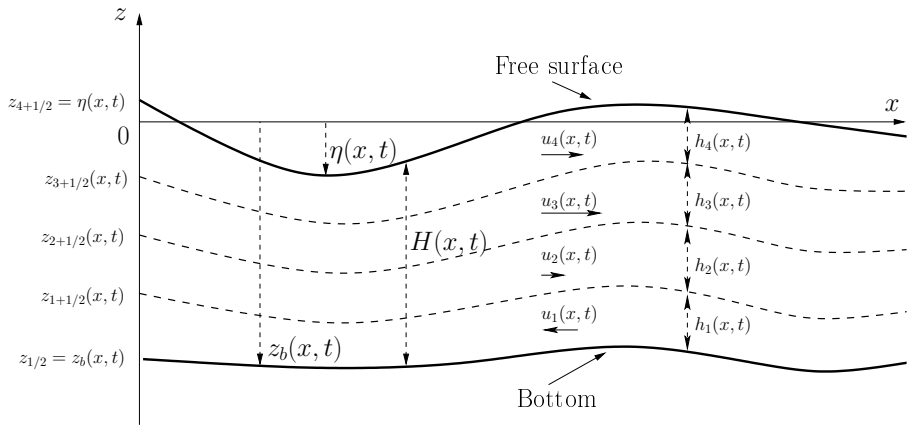

\begin{center}
\resizebox{12cm}{!}{\input Figures/notations_mc.pstex_t}
%\resizebox{12cm}{!}{\input Figures/notations_mc.pdf_t}
\caption{Notations for the multilayer approach.}
\label{fig:notations_mc}
\end{center}
\end{figure}

We depart from the free surface hydrostatic Euler system 
\begin{eqnarray}
\frac{\partial u}{\partial x }+\frac{\partial w }{\partial z }& =&0,\label{eq:div}\\
\frac{\partial u}{\partial t } + \frac{\partial u^2}{\partial x } +\frac{\partial uw }{\partial z }+ \frac{\partial p}{\partial x } &=&0,\label{eq:u1}\\
\frac{\partial p}{\partial z }& =& -g,\label{eq:p}
\end{eqnarray}
for
$$t>t_0, \quad x \in \R, \quad z_b(x) \leq z \leq \eta(x,t),$$
where $\eta(x,t)$ represents the free surface elevation, ${\bf u}=(u,w)^T$ the 
velocity. The water
height is $H = \eta - z_b$, see Fig.~\ref{fig:notations_mc}.

We add the two classical kinematic boundary conditions. 
At the free surface, we prescribe 
\begin{equation}
\frac{\partial \eta}{\partial t} + u_s \frac{\partial \eta}{\partial x}
-w_s = 0,
\label{eq:free_surf} 
\end{equation}
where the subscript $s$ denotes the value of the
considered quantity at the free surface. At the bottom, the impermeability condition gives
\begin{equation}
u_b \frac{\partial z_b}{\partial x}
-w_b = 0,
\label{eq:bottom1} 
\end{equation}
where the subscript $b$ denotes the value of the
considered quantity at the bottom.

We consider that the flow domain is divided in the vertical direction 
into $N$ layers 
of thickness $h_\alpha$ with $N+1$ interfaces $z_{\alpha+1/2}(x,t)$, $\alpha=0,...,N$
(see Fig.~\ref{fig:notations_mc}) so that  
\begin{equation}
H=\sum_{\alpha=1}^N h_\alpha,
\label{eq:Hd}
\end{equation}
and
\begin{equation}
z_{\alpha+\frac{1}{2}}(x,t) = z_b(x) + \sum_{j=1}^\alpha h_j(x,t).
\label{eq:layer}
\end{equation}

We consider the average velocities $u_\alpha$, $\alpha=1,\ldots,N$ defined by
\begin{equation}
u_\alpha(x,t) = \frac{1}{h_\alpha} \int_{z_{\alpha-1/2}}^{z_\alpha+1/2} u(x,z,t) dz,
\label{eq:u_alpha}
\end{equation}
we also denote
\begin{equation}
<u^2>_\alpha(x,t) = \frac{1}{h_\alpha} \int_{z_{\alpha-1/2}}^{z_\alpha+1/2} u^2(x,z,t) dz,
\label{eq:u2_alpha}
\end{equation}
and
\begin{equation}
u_{\alpha+1/2}=u(x,z_{\alpha+1/2},t),
\end{equation}
the value of the velocity at the interface $z_{\alpha+1/2}$. 

\begin{prpstn}
\label{prop:inviscid}
With these notations, 
an integration of (\ref{eq:div})-(\ref{eq:p}) 
over the layers $[z_{\alpha-1/2},z_{\alpha+1/2}]$, $\alpha=1,...,N$  
leads to the following system of balance laws
\begin{eqnarray}
\frac{\partial h_\alpha}{\partial t } +\displaystyle  \frac{\partial h_\alpha  u_\alpha}{\partial x}&=&G_{\alpha+1/2}- G_{\alpha-1/2},\label{eq:eq3}\\
\frac{\partial h_\alpha u_\alpha}{\partial t}
+ \frac{\partial }{\partial x}
\left(h_\alpha <u^2>_\alpha \right)+ gh_\alpha \frac{\partial H }{\partial x}  
&=&-gh_\alpha \frac{\partial  z_b}{\partial x} + u_{\alpha+1/2}G_{\alpha+1/2}-u_{\alpha-1/2}G_{\alpha-1/2}.
\label{eq:eq4}
\end{eqnarray}
The expression of the exchange terms $G_{\alpha+1/2}$ is given in the following.
\label{prop:mc_simple}
\end{prpstn}

\begin{proof}
The proof relies on simple calculus based on the Leibniz rule. 
Using the incompressibility condition (\ref{eq:div}) integrated  
over the interval $[z_{\alpha-1/2},z_{\alpha+1/2}]$, 
we deduce the mass equation (\ref{eq:eq3})
where we exhibit the kinematic of the interface on the right hand side
\begin{equation}
 G_{\alpha+1/2} = 
\frac{\partial  z_{\alpha+1/2}}{\partial t}
+ u_{\alpha+1/2} \frac{\partial  z_{\alpha+1/2}}{\partial x} 
- w(x,z_{\alpha+1/2},t),
\qquad \alpha=0,\ldots,N. \label{eq:Qalpha}
\end{equation}
The relation (\ref{eq:Qalpha}) gives the mass flux leaving/entering the layer $\alpha$ through
the interface $z_{\alpha+1/2}$.

Then we consider the velocity equation (\ref{eq:u1}).
We first observe that from the hydrostatic assumption (\ref{eq:p}) 
one can compute the pressure as a function of the water height : 
$$p(x,z,t)=g (\eta(x,t)-z).$$
Now we integrate the equation (\ref{eq:u1}) over the
interval $[z_{\alpha-1/2},z_{\alpha+1/2}]$ 
and we obtain the relation
\begin{equation}
\frac{\partial h_\alpha u_\alpha }{\partial t} + \frac{\partial }{\partial x}(h_\alpha <u^2>_\alpha) 
+ gh_\alpha \frac{\partial  \eta}{\partial x}
=u_{\alpha+1/2}G_{\alpha+1/2}-u_{\alpha-1/2}G_{\alpha-1/2},
\label{eq:mom1}
\end{equation}
and with the definition of $H$, this is equivalent to (\ref{eq:eq4}). Then the kinematic boundary conditions 
(\ref{eq:free_surf}) and (\ref{eq:bottom1}) can be written
\begin{equation}
 G_{1/2} = 0,\quad G_{N+1/2}=0.
\label{eq:Qlim}
\end{equation}
These equations just express that there is no loss/supply of mass through the bottom and the free surface.
 
Notice also that one can compute $G_{\alpha+1/2}$, just adding up the equations (\ref{eq:eq3}) for $j \leq \alpha$ and using the first equality of (\ref{eq:Qlim})
\begin{equation}
G_{\alpha+1/2}=\frac{\partial }{\partial t }\sum_{j=1}^\alpha h_j   +  \frac{\partial }{\partial x }\sum_{j=1}^\alpha h_ju_j, \qquad \alpha=1,\ldots,N. 
\label{eq:Q}
\end{equation}
\end{proof}
The standard multilayer Saint-Venant system \cite{audusse} is obtained
by prescribing
 \begin{equation}
G_{\alpha+1/2} \equiv 0. 
\label{eq:G0}
\end{equation}
This choice is clearly natural for inmiscible fluids but is not justified
if the multilayer system is seen as a numerical approximation of the hydrostaic Euler
equations. Indeed there is no reason to
prevent the water exchanges between connected layers.
Moreover it is exhibited in \cite{audusse}
that this choice may lead to the development of instabilities at the interfaces.

Here we drop this assumption and we only keep 
the two physical kinematic boundary conditions (\ref{eq:Qlim}).
 The equation~(\ref{eq:eq3}) is then no nore meaningful 
since the quantity 
$\frac{\partial h_\alpha}{\partial t},$
appears on both side of the equality. 
Nevertheless the sum of the equations~(\ref{eq:eq3})
for all the layers is still relevant
and the boundary condition (\ref{eq:Qlim}) leads to a global continuity equation
for the total water height $H$
\begin{equation}
\frac{\partial H }{\partial t} + 
\frac{\partial}{\partial x} \sum_{\alpha=1}^N h_\alpha u_\alpha = 0,
\label{eq:H}
\end{equation}
and each layer depth $h_\alpha$ is then deduced from the total water height 
by the relation 
\begin{equation}
 h_\alpha = l_\alpha H,
\label{eq:ha}
\end{equation}
with $l_\alpha$, $\alpha=1,...,N$ a given number satisfying 
\begin{equation}
l_\alpha \geq 0, \quad \sum_{\alpha=1}^N l_\alpha =1.
\label{eq:la} 
\end{equation}

Thus the momentum equation (\ref{eq:eq4}) becomes
\begin{equation}
\frac{\partial h_\alpha u_\alpha}{\partial t} 
+ \frac{\partial }{\partial x}
\left(h_\alpha <u^2>_\alpha + \frac{1}{l_\alpha} \frac{gh_\alpha^2}{2}\right) 
=-gh_\alpha \frac{\partial  z_b }{\partial x}+ u_{\alpha+1/2}G_{\alpha+1/2}-u_{\alpha-1/2}G_{\alpha-1/2}.
\label{eq:moma}
\end{equation}
Using (\ref{eq:H}),(\ref{eq:ha}), the expression of $G_{\alpha+1/2}$ given by (\ref{eq:Q}) can also be written
\begin{equation}
G_{\alpha+1/2}=\sum_{j=1}^\alpha \left( \frac{\partial h_ju_j}{\partial x }  - l_j \sum_{i=1}^N \frac{\partial h_iu_i}{\partial x } \right). 
\label{eq:Qdx}
\end{equation}
Finally we have to define the quantities $h_\alpha < u^2 >_\alpha$ and
$u_{\alpha+1/2}$ appearing in (\ref{eq:eq4}). As usual in the derivation of such systems, we have considered $h_\alpha < u^2 >_\alpha \approx h_\alpha
u_\alpha^2$, this will be discussed in details in paragraph
\ref{subsec:multilay}. The velocities $u_{\alpha+1/2}$, $\alpha=1,\ldots,N-1$ are obtained using an upwinding
\begin{equation}
u_{\alpha+1/2} =
\left\{\begin{array}{ll}
u_\alpha & \mbox{if } \;G_{\alpha+1/2} \geq 0\\
u_{\alpha+1} & \mbox{if } \;G_{\alpha+1/2} < 0.
\end{array}\right.
\label{eq:upwind}
\end{equation}
To illustrate the formulation of the new model, 
we compare it with the system proposed in \cite{audusse} 
in the simple case of a two-layer formulation.
Neglecting the viscosity and friction,
the formulation obtained by Audusse \cite{audusse} corresponds to (\ref{eq:eq3}),(\ref{eq:eq4}) with
(\ref{eq:G0}), i.e. 
\begin{eqnarray}
& & \frac{\partial h_1}{\partial t} + \frac{\partial h_1 u_1 }{\partial x} = 0,\quad\quad \frac{\partial h_2}{\partial t} + \frac{\partial h_2 u_2}{\partial x} = 0,\label{eq:audusse1}\\
& & \frac{\partial h_1 u_1}{\partial t} + \frac{\partial h_1 u_1^2}{\partial x} + gh_1 \frac{\partial (h_1+h_2)}{\partial x} = - gh_1 \frac{\partial z_b}{\partial x},\label{eq:audusse2}\\
& &\frac{\partial h_2 u_2}{\partial t} + \frac{\partial h_2 u_2^2}{\partial x} + gh_2 \frac{\partial (h_1+h_2)}{\partial x} = - gh_2 \frac{\partial z_b}{\partial x},\label{eq:audusse3}
\end{eqnarray}
with $h_1 + h_2 = H$. The preceding formulation corresponds to a superposition
of two single layer Saint-Venant systems (see also \cite{bouchut1,pares,pares1}
where a very similar model is considered in a bi-fluid framework). 

With our approach (\ref{eq:H}),(\ref{eq:moma}), the two-layer
formulation reads
\begin{eqnarray}
& & \frac{\partial H}{\partial t} + \frac{\partial h_1 u_1}{\partial x} + \frac{\partial h_2 u_2}{\partial x}  = 0,\label{eq:jsm1}\\
& & \frac{\partial h_1 u_1}{\partial t} +
\frac{\partial h_1 u_1^2}{\partial x} + \frac{g}{2}
\frac{\partial H h_1}{\partial x} = -gh_1\frac{\partial
  z_b}{\partial x} + u_{3/2}\left( l\frac{\partial H}{\partial t} + l\frac{\partial Hu_1}{\partial x} \right),\label{eq:jsm2}\\
& & \frac{\partial h_2 u_2}{\partial t} + \frac{\partial h_2 u_2^2}{\partial x} + \frac{g}{2} \frac{\partial Hh_2}{\partial x} = -gh_2\frac{\partial z_b}{\partial x} 
- u_{3/2}\left( l\frac{\partial H}{\partial t} + l\frac{\partial H
      u_1}{\partial x}\right),\label{eq:jsm3}\\
& & \mbox{where } \; h_1=l H,\quad h_2=(1-l)H, \label{eq:jsm4}
\end{eqnarray}
with $l\in(0,1)$ prescribed. 
The velocity at the interface, denoted $u_{3/2}$,
is calculated using upwinding, following the sign of the mass exchange between the layers. 
It is important to notice that, in the new formulation
(\ref{eq:jsm1})-(\ref{eq:jsm4}), we obtain directly a left hand side term written in conservative form with the topography and the mass exchange as source terms whereas the pressure term of
(\ref{eq:audusse1})-(\ref{eq:audusse3}) has to be modified \cite{audusse} to get a conservative form.
Moreover we prove in Section \ref{sec:properties} that the system (\ref{eq:jsm1})-(\ref{eq:jsm4}) 
is hyperbolic, which is not the case for system (\ref{eq:audusse1})-(\ref{eq:audusse3}).

The difference between
(\ref{eq:jsm1})-(\ref{eq:jsm4}) and (\ref{eq:audusse1})-(\ref{eq:audusse3})
mainly comes from the physical definition of the layers. Audusse
introduces a physical discretization 
where each layer has its own continuity equation. 
These $N$ continuity equations mean the layers are isolated each other, this
situation corresponds to the case of $N$ non miscible fluids. In the
formulation (\ref{eq:jsm1})-(\ref{eq:jsm4}), the discretization correponds
to a finite elements approximation --~of $P_0$ type~-- of the velocity
$u$. In this case, the definition of the layers does not correspond to a
physical partition of the flow but is related to the quality of the
desired approximation over $u$. Thus we have only one continuity
equation meaning the fluid can circulate form one layer to another.

%%%%%%%%%%%%%%%%%%%%%%%%%%%%%%%%%%%%%%%%%%%%%%%%%%%%%%%%%%%%%%%%%%%%%%%%%%%%%%%%%%%%%%%%%%%%%%%
\section{Derivation of the viscous multilayer shallow water system}
\label{sec:multi}
%%%%%%%%%%%%%%%%%%%%%%%%%%%%%%%%%%%%%%%%%%%%%%%%%%%%%%%%%%%%%%%%%%%%%%%%%%%%%%%%%%%%%%%%%%%%%%%
In this section we will apply to the Navier-Stokes equations the
multilayer approach presented in the preceding section.
%%%%%%%%%%%%%%%%%%%%%%%%%%%%%%%%%%%%%%%%%%%%%%%%%%%%%%%%%%%%%%%%%%%%%%%%%%%%%%%%%%%%%%%%%%%%%%%
\subsection{The Navier-Stokes equations}
\label{subsec:NS}
Let us start with the incompressible Navier-Stokes system \cite{lions} 
restricted to two dimensions
with gravity in which the $z$ axis represents the vertical direction. 
For simplicity, the viscosity will be kept constant and isotropic throughout the paper
(we refer the reader to \cite{decoene} for a more general framework). 
Therefore we have the following general formulation:
\begin{eqnarray}
& & \frac{\partial u}{\partial x} + \frac{\partial w}{\partial z} = 0,\label{eq:NS_2d1}\\
& & \frac{\partial u}{\partial t} + u \frac{\partial u}{\partial x} + w \frac{\partial u}{\partial z} + \frac{\partial p}{\partial x} = \frac{\partial \Sigma_{xx}}{\partial x} + \frac{\partial \Sigma_{xz}}{\partial z},\label{eq:NS_2d2}\\
& & \frac{\partial w}{\partial t} + u\frac{\partial w}{\partial x} + w\frac{\partial w}{\partial z} + \frac{\partial p}{\partial z} = -g + \frac{\partial \Sigma_{zx}}{\partial x} + \frac{\partial \Sigma_{zz}}{\partial z},
\label{eq:NS_2d3}
\end{eqnarray}
and we consider this system for
$$t>t_0, \quad x \in \R, \quad z_b(x,t) \leq z \leq \eta(x,t).$$
We use the same notations as in the previous section.
We now consider the
bathymetry $z_b$ can vary with respect to abscissa $x$ and also with
respect to time $t$. The chosen form of the viscosity tensor is symetric
\begin{eqnarray*}
 \Sigma_{xx} = 2 \nu \frac{\partial u}{\partial x}, & & \Sigma_{xz} = \nu \bigl( \frac{\partial u}{\partial z} + \frac{\partial w}{\partial x} \bigr),\label{eq:visco1}\\
\Sigma_{zz} = 2 \nu \frac{\partial w}{\partial z}, && \Sigma_{zx} = \nu \bigl( \frac{\partial u}{\partial z} + \frac{\partial w}{\partial x}\bigr),
\label{eq:visco2}
\end{eqnarray*}
with $\nu$ the viscosity coefficient.

%%%%%%%%%%%%%%%%%%%%%%%%%%%%%%%%%%%%%%%%%%%%%%%%%%%%%%%%%%%%%%%%%%%%%%%%%%%%%%%%%%%%%%%%%%%%%%%
\subsection{Boundary conditions}

The system (\ref{eq:NS_2d1})-(\ref{eq:NS_2d3}) is complete with boundary conditions. The
outward and upward unit normals to the free surface ${\bf n}_s$ and to the bottom ${\bf n}_b$ are
given by
$$
{\bf n}_s = \frac{1}{\sqrt{1 + \bigl(\frac{\partial \eta}{\partial x}\bigr)^2}}  
\left(\begin{array}{c} - {\frac{\partial \eta}{\partial x}}\\ {1} \end{array} \right), 
\quad {\bf n}_b = \frac{1}{\sqrt{1 + \bigl(\frac{\partial z_b}{\partial x}\bigr)^2}}  
\left(\begin{array}{c} - {\frac{\partial z_b}{\partial x}}\\ {1} \end{array} \right).$$
Let
$\Sigma_T$ be the total stress tensor with
$$\Sigma_T = - p I_d + \left(\begin{array}{cc} \Sigma_{xx} & \Sigma_{xz}\\ \Sigma_{zx} & \Sigma_{zz}\end{array}\right).$$

At the free surface we have the kinematic boundary condition
(\ref{eq:free_surf}). Considering the air viscosity is
negligible, the continuity
of stresses at the free boundary imposes
\begin{equation}
\Sigma_T {\bf n}_s = -p^a {\bf n}_s,
\label{eq:BC_h}
\end{equation}
where $p^a=p^a(x,t)$ is a given function corresponding to the
atmospheric pressure. Relation~(\ref{eq:BC_h}) is
equivalent to
$${\bf n}_s . \Sigma_T {\bf n}_s = -p^a,\quad \hbox{and}\quad {\bf t}_s  .\Sigma_T {\bf n}_s = 0,$$
${\bf t}_s$ being orthogonal to ${\bf n}_s$.

Since we now consider the bottom can vary with respect to time $t$, the
kinematic boundary condition reads
\begin{equation}
\frac{\partial z_b}{\partial t} + u_b \frac{\partial z_b}{\partial x}
-w_b = 0,
\label{eq:bottom} 
\end{equation}
where 
$(x,t) \mapsto z_b(x,t)$ is a
given function. Notice that the equation (\ref{eq:bottom})
reduces to a classical no-penetration condition (\ref{eq:bottom1}) when $z_b$
does not depend on time $t$. For the stresses at the bottom we consider a wall law under the form
\begin{equation}
\Sigma_T {\bf n}_b - ({\bf n}_b . \Sigma_T {\bf n}_b){\bf n}_b = \kappa({\bf v_b},H) {\bf v}_b,
\label{eq:BC_z_b}
\end{equation}
with ${\bf v}_b={\bf u}_b - (0,\frac{\partial z_b}{\partial t})^T$ the
relative velocity between the water and the bottom. If $\kappa({\bf v_b},H)$ is constant then
we recover a Navier friction condition as in \cite{gerbeau}. Introducing
$k_l$ laminar and $k_t$ turbulent friction coefficients, we use the expression 
$$\kappa({\bf v_b},H) = k_l + k_t H |{\bf v_b}|,$$
corresponding to the boundary condition used in \cite{marche}. Another
form of $\kappa({\bf v_b},H)$ is used in \cite{bouchut} and for other
wall laws, the reader can also refer
to \cite{valentin}. Due to thermomechanical considerations, in the
sequel we suppose $\kappa({\bf v_b},H) \geq 0$ and $\kappa({\bf
v_b},H)$ is often simply denoted by $\kappa$.

Let ${\bf t}_b$ satisfying ${\bf t}_b . {\bf n}_b = 0$
then when multiplied by ${\bf t}_b$ and ${\bf n}_b$, Equation (\ref{eq:BC_z_b}) leads to
$${\bf t}_b. \Sigma_T {\bf n}_b = \kappa {\bf v}_b . {\bf t}_b,\quad \hbox{and}\quad {\bf v}_b. {\bf n}_b = 0.$$
%%%%%%%%%%%%%%%%%%%%%%%%%%%%%%%%%%%%%%%%%%%%%%%%%%%%%%%%%%%%%%%%%%%%%%%%%%%%%%%%%%%%%%%%%%%%%%%
\subsection{The rescaled system}
\label{subsec:NS_rescaled}
The physical system is rescaled using the quantities
\begin{itemize}
\item $h$ and $\lambda$, two characteristic dimensions along the $z$ and $x$
      axis respectively,
\item $a_s$ the typical wave amplitude, $a_b$ the typical bathymetry variation,
\item $C=\sqrt{gh}$ the typical horizontal wave speed.
\end{itemize}
Classically for the derivation of the Saint-Venant system, we introduce
the small parameter
$$\varepsilon = \frac{h}{\lambda}.$$
When considering long waves propagation,
another important parameter needs be considered, namely
$$\delta = \frac{a_s}{h},$$
and we consider for the bathymetry $\frac{a_b}{h} = {\cal
O}(\delta)$. Notice that $\varepsilon$ is related to {\it a priori}
informations only associated to geometrical features whereas $a_s$ and
accordingly $\delta$ deal with the state variables of the problem.

Depending on the application, $\delta$ can be considered or not as a small
parameter. For finite amplitude wave theory and assuming
$z_b(x,t)=z_b^0$, one considers $\varepsilon \ll 1$, $\delta =
{\cal O}(1)$ whereas the Boussinesq waves theory requires
$$\delta \ll 1,\quad\varepsilon \ll 1\quad\hbox{and}\quad U_r = {\cal O}(1)$$
where $U_r$ is the Ursell number defined by $U_r = \frac{\delta}{\varepsilon^2}$, see \cite{ursell}. All along this work, we consider $\varepsilon \ll 1$ whereas, even if the parameter $\delta$ is introduced in the
rescaling, the assumption $\delta \ll 1$ is not considered except when
explictly mentioned.

As for the Saint-Venant system
\cite{gerbeau,marche}, we introduce some characteristic quantities~:
$T=\lambda/C$ for the time, $W = a_s/T 
= \varepsilon\delta C$
for the vertical velocity, $U = W/\varepsilon=\delta C$, for the horizontal velocity, $P= C^2$ for
the pressure. This leads to the following dimensionless quantities
$$\tilde{x} = \frac{x}{\lambda},\quad \tilde{z} = \frac{z}{h},\quad\tilde{\eta} = \frac{\eta}{a_s},\quad\tilde{t} = \frac{t}{T},$$
$$\tilde{p} = \frac{p}{P},\quad \tilde{u} = \frac{u}{
U},\quad\mbox{and}\quad \tilde{w} = \frac{w}{W}.$$
Notice that the definition of the charateristic velocities implies $\delta = \frac{U}{C}$ so $\delta$ also corresponds to the Froude number. When $\delta = {\cal O}(1)$ we have $U \approx C$ and we recover the classical rescaling used for the Saint-Venant system. For the bathymetry $z_b$ we write $z_b(x,t) = Z_b(x)+b(t)$ and we
introduce $\tilde{z}_b = Z_b/h$ and $\tilde{b} = b / a_b$. This leads to
$$\frac{\partial z_b}{\partial t} = \varepsilon\delta C
\frac{\partial\tilde{b}}{\partial \tilde{t}} = W\frac{\partial\tilde{b}}{\partial \tilde{t}},\quad\mbox{and}\quad\frac{\partial z_b}{\partial x} = \varepsilon\frac{\partial\tilde{z}_b}{\partial \tilde{x}}.$$
The different rescaling applied to the time
and space derivatives of $z_b$ means that a classical shallow water
assumption is made concerning the space variations of the bottom profile whereas we assume
the time variations of $z_b$ lie in the framework of the Boussinesq
assumption and are consistent with the rescaling applied to the velocity
$w$.

We also introduce $\tilde{\nu} = \frac{\nu}{\lambda C}$ and we set $\tilde{\kappa} = \frac{\kappa}{C}$.
Notice that the definitions for the dimensionless quantities are
consistent with the one used for the Boussinesq
system \cite{peregrine,walkley}. Notice also that the rescaling used by
Nwogu \cite{nwogu} differs from the preceding one since Nwogu uses
$\tilde{w} = \frac{\varepsilon^2}{W}w$.

As in \cite{gerbeau,marche}, we suppose we are in the following asymptotic regime
$$\tilde{\nu} = \varepsilon \nu_0,\qquad \mbox{and}\qquad\tilde{\kappa} =
\varepsilon \kappa_0,$$
with $\kappa_0 = \kappa_{l,0} + \varepsilon\kappa_{t,0}(\tilde{\bf
v}_b,\tilde{H})$, $\kappa_{l,0}$ being constant.

This non-dimensionalization of the Navier-Stokes system (\ref{eq:NS_2d1})-(\ref{eq:NS_2d3}) leads to
\begin{eqnarray}
&&\frac{\partial \tilde{u}}{\partial \tilde{x}} + \frac{\partial \tilde{w}}{\partial \tilde{z}} = 0,\label{eq:NS_2d_newa1}\\
&&\varepsilon\delta\frac{\partial \tilde{u}}{\partial \tilde{t}} + \varepsilon\delta^2 \frac{\partial 
{\tilde{u}}^2}{\partial \tilde{x}} + \varepsilon\delta^2\frac{\partial \tilde{u} \tilde{w} }{\partial \tilde{z}} + \varepsilon\frac{\partial \tilde{p}}{\partial \tilde{x}} = \varepsilon^2\delta\frac{\partial }{\partial
\tilde{x}}\left(2\nu_0\frac{\partial \tilde{u}}{\partial \tilde{x}}\right)
 + \frac{\partial }{\partial \tilde{z}}\left( \delta\nu_0\frac{\partial \tilde{u}}{\partial \tilde{z}} + 
\varepsilon^2\delta\nu_0\frac{\partial \tilde{w}}{\partial \tilde{x}}\right), \label{eq:NS_2d_newa2}\\
&&\varepsilon^2\delta\left(\frac{\partial \tilde{w}}{\partial \tilde{t}} + \delta\frac{\partial 
\tilde{u} \tilde{w}}{\partial \tilde{x}} + \delta\frac{\partial {\tilde{w}}^2}{\partial \tilde{z}} \right) + \frac{\partial 
\tilde{p}}{\partial \tilde{z}} = -1 
 + \frac{\partial }{\partial \tilde{x}}\left( \varepsilon\delta\nu_0\frac{\partial \tilde{u}}{\partial \tilde{z}} + 
\nu_0\varepsilon^3\delta \frac{\partial \tilde{w}}{\partial \tilde{x}}\right)
 + \varepsilon\delta\frac{\partial }{\partial \tilde{z}}\left(2\nu_0\frac{\partial \tilde{w}}{\partial \tilde{z}}\right), \label{eq:NS_2d_newa3}
\end{eqnarray}
where we use the divergence free condition
to write velocity equations (\ref{eq:NS_2d_newa2}) and
(\ref{eq:NS_2d_newa3}) in a conservative form. The associated boundary conditions (\ref{eq:free_surf}), (\ref{eq:BC_h}), 
(\ref{eq:bottom}) and (\ref{eq:BC_z_b}) become
\begin{eqnarray}
\lefteqn{\frac{\partial \tilde{\eta}}{\partial \tilde{t}} + \delta\tilde{u}_s \frac{\partial \tilde{\eta}}{\partial \tilde{x}} - \tilde{w}_s = 0,\label{eq:BC_re1}}\\
\lefteqn{2\varepsilon\delta\nu_0\left.\frac{\partial \tilde{w}}{\partial \tilde{z}}\right|_s -\tilde{p}_s -\varepsilon\delta^2\nu_0\frac{\partial \tilde{\eta}}{\partial \tilde{x}}\left(\left.\frac{\partial \tilde{u}}{\partial \tilde{z}}\right|_s + \varepsilon^2\left.\frac{\partial \tilde{w}}{\partial \tilde{x}}\right|_s\right) = -\delta\tilde{p}^a,\label{eq:BC_re2}}\\
\lefteqn{\delta\nu_0\left(\left.\frac{\partial \tilde{u}}{\partial \tilde{z}}\right|_s + \varepsilon^2\left.\frac{\partial \tilde{w}}{\partial \tilde{x}}\right|_s \right) - \varepsilon\delta\frac{\partial \tilde{\eta}}{\partial \tilde{x}} \left(2\varepsilon\delta\nu_0\left.\frac{\partial \tilde{u}}{\partial \tilde{x}}\right|_s - \tilde{p}_s\right) = \varepsilon\delta^2\frac{\partial \tilde{\eta}}{\partial \tilde{x}} \tilde{p}^a, \label{eq:BC_re3}}\\
\lefteqn{\frac{\partial \tilde{b}}{\partial \tilde{t}} + \tilde{u}_b \frac{\partial \tilde{z}_b}{\partial \tilde{x}} - \tilde{w}_b = 0, \label{eq:BC_re4}}\\
\lefteqn{\delta\nu_0\left(\varepsilon^2\left.\frac{\partial \tilde{w}}{\partial \tilde{x}}\right|_b + \left.\frac{\partial \tilde{u}}{\partial \tilde{z}}\right|_b \right) - \varepsilon\frac{\partial \tilde{z}_b}{\partial \tilde{x}} \left(2\varepsilon\delta\nu_0 \left.\frac{\partial \tilde{u}}{\partial \tilde{x}}\right|_b - p_b \right) \nonumber}\\
& & \quad + \varepsilon\frac{\partial \tilde{z}_b}{\partial \tilde{x}}\left(2\varepsilon\delta\nu_0\left.\frac{\partial \tilde{w}}{\partial \tilde{z}}\right|_b\right. - p_b - \left.\varepsilon\nu_0\frac{\partial \tilde{z}_b}{\partial \tilde{x}} \left(\delta\left.\frac{\partial \tilde{u}}{\partial \tilde{z}}\right|_b + \varepsilon^2\delta\left.\frac{\partial \tilde{w}}{\partial \tilde{x}}\right|_b\right)\right)\nonumber\\
& & = \varepsilon\delta\kappa_0 \sqrt{1 + \varepsilon^2\left(\frac{\partial \tilde{z}_b}{\partial \tilde{x}}\right)^2} \left(\tilde{u}_b + \varepsilon^2\frac{\partial \tilde{z}_b}{\partial \tilde{x}}\bigl(\tilde{w}_b - \frac{\partial \tilde{b}}{\partial \tilde{t}}\bigr)\right) \label{eq:BC_re5}.
\end{eqnarray}
For the sake of clarity, in the sequel we drop the symbol $\tilde{}$ and
we denote $\frac{\partial b}{\partial t} = \frac{\partial z_b}{\partial t}$.

%%%%%%%%%%%%%%%%%%%%%%%%%%%%%%%%%%%%%%%%%%%%%%%%%%%%%%%%%%%%%%%%%%%%%%%%%%%%%%%%%%%%%%%%%%%%%%%
\subsection{The Shallow Water system}
\label{sec:SW}

The derivation of multilayer approximation is somehow technical. In order to better explain the analysis we recall the monolayer case following the asymptotic expansion in \cite{gerbeau}.

In the following the two sets of equations
(\ref{eq:NS_2d_newa1})-(\ref{eq:NS_2d_newa3}) and (\ref{eq:BC_re1})-(\ref{eq:BC_re5}) are approximated to retain only
the high order terms. 

Due to the applied rescaling some terms of the viscosity tensor e.g.
$$\varepsilon^3\delta\frac{\partial }{\partial x}\left( \nu_0 \frac{\partial w}{\partial x}\right)$$
are very small and could be neglected. But, as mentioned in
\cite[Remarks~1 and 2]{audusse}, the approximation of the viscous terms
has to preserve the dissipation energy that is an essential property of
the Navier-Stokes and averaged Navier-Stokes equations. Since we
privilege this stability requirement and in order to keep a symmetric
form of the viscosity tensor, we consider in the sequel a modified version of (\ref{eq:NS_2d_newa1})-(\ref{eq:NS_2d_newa3}) under the form
\begin{eqnarray}
&&\frac{\partial u}{\partial x} + \frac{\partial w}{\partial z} = 0, \label{eq:NS_2d_new1}\\
&&\varepsilon\delta\frac{\partial u}{\partial t} + \varepsilon\delta^2 \frac{\partial u^2}{\partial x} + \varepsilon\delta^2 \frac{\partial u w}{\partial z} + \varepsilon\frac{\partial p}{\partial x} = \varepsilon^2\delta\frac{\partial }{\partial x}\left(2\nu_0\frac{\partial u}{\partial x}\right)
+ \frac{\partial }{\partial z}\left( \delta\nu_0\frac{\partial u}{\partial z}\right),\label{eq:NS_2d_new2}\\
&&\varepsilon^2\delta\left(\frac{\partial w}{\partial t} + \delta \frac{\partial u w}{\partial x} + \delta\frac{\partial w^2}{\partial z} \right) + \frac{\partial p}{\partial z} = -1 + \frac{\partial }{\partial x}\left( \varepsilon\delta\nu_0\frac{\partial u}{\partial z}\right)
+ \frac{\partial }{\partial z}\left( 2\varepsilon\delta\nu_0\frac{\partial w}{\partial z}\right),\label{eq:NS_2d_new3}
\end{eqnarray}
corresponding to a viscosity tensor of the form
$$\Sigma_{xx} = 2 \nu \frac{\partial u}{\partial x}, \quad \Sigma_{xz} =
\Sigma_{zx} = \nu\frac{\partial u}{\partial z}, \quad \Sigma_{zz} = 2
\nu \frac{\partial w}{\partial z}.$$
This means the terms in
$\varepsilon^2  {\partial_x w}$
have been neglected in (\ref{eq:NS_2d_newa1})-(\ref{eq:NS_2d_newa3}) and in
(\ref{eq:BC_re1})-(\ref{eq:BC_re5}). For details about the adopted form of the viscosity tensor see \cite[Remark~2]{JSM_DCDS} and \cite[Lemma~2.1]{audusse}.

In the same way, retaining only the high order terms, the boundary conditions (\ref{eq:BC_re1})-(\ref{eq:BC_re5})
become
\begin{eqnarray}
&&\frac{\partial {\eta}}{\partial {t}} + \delta{u}_s \frac{\partial {\eta}}{\partial {x}} - {w}_s = 0,\label{eq:BC_ap1}\\
&&2\varepsilon\delta\nu_0\left.\frac{\partial {w}}{\partial {z}}\right|_s -{p}_s -\varepsilon\delta^2\nu_0\frac{\partial {\eta}}{\partial {x}}\left.\frac{\partial {u}}{\partial {z}}\right|_s = -\delta{p}^a,\label{eq:BC_ap2}\\
&&\delta\nu_0\left.\frac{\partial {u}}{\partial {z}}\right|_s - \varepsilon\delta\frac{\partial {\eta}}{\partial {x}} \left(2\varepsilon\delta\nu_0\left.\frac{\partial {u}}{\partial {x}}\right|_s - {p}_s\right) = \varepsilon\delta^2\frac{\partial {\eta}}{\partial {x}} {p}^a, \label{eq:BC_ap3}\\
&&\frac{\partial {z_b}}{\partial {t}} + {u}_b \frac{\partial {z}_b}{\partial {x}} - {w}_b = 0, \label{eq:BC_ap4}\\
&&\delta\nu_0\left.\frac{\partial {u}}{\partial {z}}\right|_b - \varepsilon\frac{\partial {z}_b}{\partial {x}} \left(2\varepsilon\delta\nu_0 \left.\frac{\partial {u}}{\partial {x}}\right|_b - p_b \right) 
+ \varepsilon\frac{\partial {z}_b}{\partial {x}}\left(2\varepsilon\delta\nu_0\left.\frac{\partial {w}}{\partial {z}}\right|_b\right. - p_b - \left.\varepsilon\delta\nu_0\frac{\partial {z}_b}{\partial {x}}\left.\frac{\partial {u}}{\partial {z}}\right|_b\right)\nonumber\\
&&\qquad = \varepsilon\delta\kappa_0 \left(1 + \varepsilon^2\left(\frac{\partial {z}_b}{\partial {x}}\right)^2\right)^{3/2} u_b \label{eq:BC_ap5}.
\end{eqnarray}
Now we will exhibit the hydrostatic and non hydrostatic parts of the pressure.
An integration of (\ref{eq:NS_2d_new3}) from $z$ to $\delta \eta$ gives
\begin{eqnarray}
&&\varepsilon^2\delta \int_z^{\delta \eta} \bigl(\frac{\partial w}{\partial t} + \delta\frac{\partial (uw)}{\partial x}\bigr) dz + \varepsilon^2\delta^2(w_s^2 - w^2) +
p_s - p \nonumber\\
&&\qquad = -(\delta \eta - z)+ \varepsilon\delta\int_z^{\delta \eta} \frac{\partial}{\partial x}\left(\nu_0 \frac{\partial u}{\partial
z}\right)dz - 2\varepsilon\delta\nu_0\frac{\partial w}{\partial z} + \left.2\varepsilon\delta\nu_0\frac{\partial w}{\partial z}\right|_s. \label{eq:p_s}
\end{eqnarray}
From the equations (\ref{eq:BC_re2}) and (\ref{eq:BC_re3}) it comes
\begin{equation}
\left.\frac{\partial u}{\partial z}\right|_s = {\cal O}(\varepsilon^2),
\label{eq:bcdus}
\end{equation}
and the boundary condition (\ref{eq:BC_re2}) gives
\begin{equation}
p_s = \delta p^a + 2\varepsilon\delta\left.\frac{\partial w}{\partial
z}\right|_s + {\cal O}(\varepsilon^3\delta^2).
\label{eq:bcp_s}
\end{equation}
The previous relation and the kinematic boundary condition
(\ref{eq:BC_re1}) allow us to rewrite (\ref{eq:p_s}) under the form
\begin{eqnarray*}
&&\varepsilon^2\delta  \left(\frac{\partial}{\partial t}\int_z^{\delta \eta} w\ dz + \delta\frac{\partial }{\partial x}\int_z^{\delta \eta} (uw)\ dz \right) - \varepsilon^2\delta^2 w^2 +
\delta p^a - p \nonumber\\
 &&\qquad = -(\delta \eta - z)+ \varepsilon\delta\int_z^{\delta \eta} \frac{\partial}{\partial x}\left(\nu_0 \frac{\partial u}{\partial
z}\right)dz - 2\varepsilon\delta\nu_0\frac{\partial w}{\partial z} + {\cal O}(\varepsilon^3\delta).
\label{eq:NS3_1}
\end{eqnarray*}
Classically we have
\begin{equation}
\frac{\partial u_s}{\partial x} = \left.\frac{\partial u}{\partial
x}\right|_s + \delta\frac{\partial \eta}{\partial
x}\left.\frac{\partial u}{\partial z}\right|_s = \left.\frac{\partial
u}{\partial x}\right|_s + {\cal O}(\varepsilon^2\delta),
\label{eq:deriv_comp}
\end{equation}
and using relations (\ref{eq:NS_2d_new1}), (\ref{eq:deriv_comp}) and the
Leibniz rule we have
\begin{equation*}
\varepsilon\delta\int_z^{\delta \eta} \frac{\partial}{\partial x}\left(\nu_0 \frac{\partial u}{\partial
z}\right) dz - 2\varepsilon\delta\nu_0\frac{\partial w}{\partial z} 
= \varepsilon\delta\nu_0\frac{\partial u}{\partial x} + \left.\varepsilon\delta\nu_0\frac{\partial u}{\partial x}\right|_s + {\cal O}(\varepsilon^3\delta).
\end{equation*}
This leads to the expression for the pressure $p$
\begin{eqnarray}
p = p_{h} + p_{nh} + {\cal O}(\varepsilon^3\delta),
\label{eq:pfin}
\end{eqnarray}
where the viscous and hydrostatic part $p_{h}$ is given by
$$p_h = \delta p^a + (\delta \eta - z) -
\varepsilon\delta\nu_0\frac{\partial u}{\partial x} -
\varepsilon\delta\nu_0\left.\frac{\partial u}{\partial x}\right|_s,$$
and the non-hydrostatic part $p_{nh}$ is
$$ p_{nh} = \varepsilon^2\delta  \left(\frac{\partial}{\partial
t}\int_z^{\delta \eta} w\ dz + \delta\frac{\partial }{\partial
x}\int_z^{\delta \eta} (uw)\ dz \right) - \varepsilon^2\delta^2 w^2.$$
The derivation and analysis of a classical Saint-Venant type system taking into account
the non-hydrostatic part of the pressure has already been carried out by the authors \cite{JSM_DCDS}. 
The derivation of the multilayer system in this general framework is in progress. 
It will be presented in a forthcoming paper.

In the sequel, we restrict to the situation $p_{nh}=0$.
Due to this hydrostatic assumption, we have 
\begin{equation}
p=p_h+{\cal O}(\varepsilon^2\delta),
\label{eq:p_heps}
\end{equation}
and we retain for $p_h$ the expression
\begin{equation}
p_h = \delta p^a +(\delta \eta - z) -
2\varepsilon\delta\nu_0\frac{\partial u}{\partial x}.
\label{eq:p_h}
\end{equation}
%%%%%%%%%%%%%%%%%%%%%%%%%%%%%%%%%%%%%%%%%%%%%%%%%%%%%%%%%%%%%%%%%%%%%%%%%%%%%%%%%%%%%%%%%%%%%%%
Then using (\ref{eq:BC_re3}), (\ref{eq:BC_re5}) and (\ref{eq:bcp_s}) one obtains
\begin{equation}
\left.\frac{\partial u}{\partial z}\right|_s = {\cal O}(\varepsilon^2),\quad \left.\frac{\partial u}{\partial z}\right|_b = {\cal O}(\varepsilon).
\label{eq:BC1_ubis}
\end{equation}
From (\ref{eq:p_heps}),(\ref{eq:p_h}) we can write
\begin{equation}
p - \delta p^a = \delta\eta - z + {\cal O}(\varepsilon\delta),
\label{eq:p_simple}
\end{equation}
leading to
$$\frac{\partial p}{\partial x} = {\cal O}(\delta).$$
The preceding relation inserted in (\ref{eq:NS_2d_new2}) leads to
\begin{equation}
\nu_0\frac{\partial^2 u}{\partial z^2} = {\cal O}(\varepsilon),
\label{eq:BC1_uter}
\end{equation}
and Equations (\ref{eq:BC1_ubis}) and (\ref{eq:BC1_uter}) mean that
\begin{equation}
u(x,z,t) = u(x,0,t) + {\cal O}(\varepsilon),
\label{eq:u}
\end{equation}
i.e. we recognize the so-called ``motion by slices'' of the usual Saint-Venant system. 
If we introduce the averaged quantity
$$\bar{u} = \frac{1}{\delta \eta - z_b} \int_{z_b}^{\delta \eta} u\
dz,$$
it is well known \cite{JSM_DCDS,saleri,gerbeau,marche} that the shallow water system (\ref{eq:NS_2d_new1}),(\ref{eq:NS_2d_new2}) 
with an hydrostatic pressure (\ref{eq:p_heps}),(\ref{eq:p_h}) is approximated in ${\cal O}(\varepsilon^2 \delta)$ by the following Saint-Venant system written with the variables with dimension
\begin{eqnarray}
&&\frac{\partial H}{\partial t} + \frac{\partial H\bar{u}}{\partial x}  = 0,\label{eq:sv1}\\
&&\frac{\partial H\bar{u}}{\partial t} + \frac{\partial H\bar{u}^2}{\partial x} + \frac{g}{2}\frac{\partial H^2}{\partial x} = -H \frac{\partial p^a}{\partial x} -gH \frac{\partial z_b}{\partial x} + \frac{\partial}{\partial x}\bigl(4 \nu H\frac{\partial \bar{u}}{\partial x}\bigr)
- \frac{\kappa(\bar{{\bf v}},H)}{1 + \frac{\kappa(\bar{{\bf v}},H)}{3\nu}H} \bar{u}\hspace*{3cm}\label{eq:sv2}
\end{eqnarray}
with $ H= \eta-z_b$.
%%%%%%%%%%%%%%%%%%%%%%%%%%%%%%%%%%%%%%%%%%%%%%%%%%%%%%%%%%%%%%%%%%%%%%%%%%%%%%%%%%%%%%%%%%%%%%%
\subsection{The viscous multilayer Shallow Water system}
\label{subsec:multilay}

We again consider the Shallow Water system (\ref{eq:NS_2d_new1}),(\ref{eq:NS_2d_new2}) with an hydrostatic pressure (\ref{eq:p_heps}),(\ref{eq:p_h}).
Here another approximation is introduced concerning the velocity $u$, it is no more assumed constant along the vertical but
is discretized in the $z$ direction using piecewise
constant functions, see Fig.~\ref{fig:notations_mc}. As introduced in Sec.~\ref{sec:simplified} the interval $[z_b,\delta\eta]$ is
divided into $N$ layers 
of thickness $h_\alpha$ and we use the definitions
(\ref{eq:ha}),(\ref{eq:la}). We write
\begin{equation}
u^{mc}(x,z,\{z_\alpha\},t)  = \sum_{\alpha=1}^N 1_{[z_{\alpha-1/2},z_{\alpha+1/2}]}(z)u_\alpha(x,t)
\label{eq:ulayer}
\end{equation}
with the velocities $u_\alpha$, $\alpha \in [1,\ldots,N]$ defined by
(\ref{eq:u_alpha}).

Notice that from (\ref{eq:layer}) we have $z_{1/2}=z_b={\cal O}(1)$ and
$z_{N+1/2}=\delta\eta={\cal O}(\delta)$. The difference of magnitude between
$z_{1/2}$ and $z_{N+1/2}$ makes the assumption $\delta \ll 1$
difficult to integrate in the definition of the
$\{z_{\alpha +1/2}\}$.

Now we try to quantify the error between $u$ and its piecewise
approximation $u^{mc}$. First we notice that in absence of friction at the
bottom and due to the Shallow Water
assumption, the relations (\ref{eq:BC1_ubis}) become
\begin{equation}
\left.\frac{\partial u}{\partial z}\right|_s = \left.\frac{\partial u}{\partial z}\right|_b = {\cal O}(\varepsilon^2).
\label{eq:BC1_mod}
\end{equation}
This means we can consider that except for the bottom layer, each layer
inherits the approximation (\ref{eq:BC1_mod}) i.e.
\begin{equation*}
\frac{\partial u}{\partial z} = {\cal O}(\varepsilon^2) \qquad\mbox{for } z\geq z_{3/2},
\end{equation*}
and therefore for all $\alpha>1$
\begin{equation}
u(x,z,t) - u_\alpha(x,t) = {\cal O}(\varepsilon^2),
\label{eq:u_approxi}
\end{equation}
or equivalently
\begin{equation*}
u(x,z,t) - u^{mc}(x,z,\{z_\alpha\},t) = {\cal O}(\varepsilon^2), \qquad\mbox{for } z\geq z_{3/2}.
\end{equation*}
In the bottom layer we only have
\begin{equation*}
u(x,z,t) - u_1(x,t) = {\cal O}(\varepsilon),
\end{equation*}
but as in \cite{JSM_DCDS,gerbeau}, it can be proved that we have an
approximation of the velocity through a parabolic correction
\begin{equation}
u = \left(1 + \frac{\varepsilon\kappa_0}{\nu_0} \bigl( z - z_b -
\frac{(z -z_b)^2}{2H} - \frac{H}{3}\bigr)\right) u_1 + {\cal O}(\varepsilon^2),
\label{eq:u_approx}
\end{equation}
for $z \in [z_{1/2},z_{3/2}]$. Using the discretization (\ref{eq:layer}),(\ref{eq:u_alpha}) and (\ref{eq:ulayer}) we claim 
\begin{prpstn}
%\begin{proposition}
The multilayer formulation of the Saint-Venant system defined by
\begin{eqnarray}
&&\frac{\partial H}{\partial t} + \sum_{\alpha=1}^N \frac{\partial h_\alpha u_\alpha}{\partial x}  = 0,\label{eq:sv_mc1}\\
\nonumber\\
&&\frac{\partial h_1 u_1}{\partial t} + \frac{\partial h_1 u_1^2}{\partial x} + \frac{g}{2 l_1} \frac{\partial h_1^2}{\partial x} = -h_1\frac{\partial p^a}{\partial x} -gh_1\frac{\partial z_b}{\partial x}
 + u_{3/2}G_{3/2}\nonumber\\
%\left(\frac{\partial h_{1}}{\partial t} + \frac{\partial h_1 u_1}{\partial x}\right)
&&\qquad + \frac{\partial }{\partial x}\left(4\nu h_1\frac{\partial u_1}{\partial x} \right) - 4\nu \frac{\partial z_{3/2}}{\partial x}\frac{\partial u_{3/2}}{\partial x}
+ 2\nu\frac{u_2 - u_1}{h_2+h_1} - \kappa(\bar{\bf v},H) u_1,\label{eq:sv_mc2}\\
\nonumber\\
&&\frac{\partial h_\alpha u_\alpha}{\partial t} + \frac{\partial h_\alpha u_\alpha^2}{\partial x} + \frac{g}{2l_\alpha} \frac{\partial h_\alpha^2}{\partial x} = -h_\alpha\frac{\partial p^a}{\partial x} -gh_\alpha\frac{\partial z_b}{\partial x}
 +u_{\alpha+1/2}G_{\alpha+1/2} -u_{\alpha-1/2}G_{\alpha-1/2}
\nonumber\\
&& \qquad + \frac{\partial }{\partial x}\left(4\nu h_\alpha\frac{\partial u_\alpha}{\partial x} \right) - 4\nu \left[\frac{\partial z_j}{\partial x}\frac{\partial u_j}{\partial x}\right]_{j=\alpha-1/2}^{j=\alpha+1/2}
%\nonumber\\
+ 2\nu\frac{u_{\alpha+1} - u_\alpha }{h_{\alpha+1}+h_\alpha} - 2\nu\frac{u_{\alpha} - u_{\alpha-1}}{h_{\alpha}+h_{\alpha-1}},\label{eq:sv_mc3}\\
&&\qquad \qquad \qquad \qquad \qquad \qquad \qquad \qquad \qquad \qquad \mbox{for } \alpha \in \{2,\ldots,N-1\}\nonumber\\
\nonumber\\
&&\frac{\partial h_N u_N}{\partial t} + \frac{\partial h_N u_N^2}{\partial x} + \frac{g}{2l_N} \frac{\partial h_N^2}{\partial x} = -h_N\frac{\partial p^a}{\partial x} - gh_N\frac{\partial z_b}{\partial x}
-u_{N-1/2}G_{N-1/2}\nonumber\\
&&\qquad  + \frac{\partial }{\partial x}\left(4\nu h_N\frac{\partial u_N}{\partial x} \right) + 4\nu\frac{\partial z_{N-1/2}}{\partial x}\frac{\partial u_{N-1/2}}{\partial x}
- 2\nu\frac{u_N - u_{N-1}}{h_N+h_{N-1}},\label{eq:sv_mc4}
\end{eqnarray}
with $h_\alpha=l_\alpha H(x,t)$ and $G_{\alpha+1/2}$ given by (\ref{eq:Q}),
results from a formal asymptotic approximation in
 ${\cal O}(\varepsilon^2\delta)$ coupled with a vertical discretization
 of the Navier-Stokes equations (\ref{eq:NS_2d_newa1})-(\ref{eq:NS_2d_newa3})
 with hydrostatic pressure.
\label{prop:sv_mc}
\end{prpstn}

\begin{proof}
The integration of the divergence equation (\ref{eq:NS_2d_new1}) on each layer 
has been already performed in the proof of Proposition \ref{prop:inviscid}.
We recall that the deduced layer mass equations (\ref{eq:eq3}) are not meaningful
if no hypothesis is made concerning the mass exchange term $G_{\alpha+1/2}$ defined by (\ref{eq:Q}).
We thus consider a global mass equation (\ref{eq:sv_mc1}) by adding them up.
We can also directly integrate the divergence equation from bottom to free surface 
in order to obtain equation (\ref{eq:sv_mc1}).

We now consider the horizontal velocity equation (\ref{eq:NS_2d_new2}) integrated over the
 interval $[z_{\alpha-1/2},z_{\alpha+1/2}]$. Using for each layer an approximation similar to
(\ref{eq:u_approxi}),(\ref{eq:u_approx}), we prove that
$$ \frac{1}{h_\alpha} \int_{z_{\alpha-1/2}}^{z_\alpha+1/2} u^2(x,z,t) dz = 
u_\alpha^2 +{\cal O}(\varepsilon^2).$$
In the context of the hydrostatic approximation, we assume that the pressure satisfies 
(\ref{eq:p_heps}),(\ref{eq:p_h}). The treatment of the inviscid part of the pressure has already been presented 
in the proof of Proposition \ref{prop:inviscid} where we have written
 for the gravitational part of the pressure
\begin{equation}
\begin{split}
\int_{z_{\alpha-1/2}}^{z_\alpha+1/2}  \frac{\partial }{\partial x}
 (\delta \eta - z)\ dz = \frac{1}{2l_\alpha}\frac{\partial }{\partial x}
 h_\alpha^2 + \frac{h_\alpha}{l_\alpha}\frac{\partial z_b}{\partial x}.
\end{split}
\label{eq:p1}
\end{equation}
Notice that it is also possible to write 
\begin{equation}
\begin{split}
\int_{z_{\alpha-1/2}}^{z_\alpha+1/2}  \frac{\partial }{\partial x} (\delta \eta - z)\ dz  = \frac{1}{2}\frac{\partial }{\partial x} \left[ h_\alpha\left( 2\sum_{j=\alpha+1}^N h_j + h_\alpha\right) \right]
 -  \frac{\partial z_{\alpha+1/2}}{\partial x} \sum_{j=\alpha+1}^N h_j + \frac{\partial z_{\alpha-1/2}}{\partial x} \sum_{j=\alpha}^N h_j.
\end{split}
\label{eq:p2}
\end{equation}
The expressions (\ref{eq:p1}) and (\ref{eq:p2}) lead to the same
 property for the complete model even if the hyperbolic part is
 modified. The second formulation seems more adapted to the physical
 description ``by layers'' of the system but leads to complementary
 source terms whose discretization is subtle. We will use and analyse
 (\ref{eq:p2}) in a forthcoming paper. In the following we use (\ref{eq:p1}).

The integration of the viscous part of the pressure leads to
$$   \int_{z_{\alpha-1/2}}^{z_\alpha+1/2} 
\frac{\partial}{\partial x} \left(2\nu\frac{\partial u}{\partial x}\right)
= 
\frac{\partial }{\partial x}
\left(2\nu h_\alpha\frac{\partial u_\alpha}{\partial x} \right) +
2\nu 
\left[\frac{\partial z_j}{\partial x}\frac{\partial u_j}{\partial x}\right]_{j=\alpha-1/2}^{j=\alpha+1/2}
+ {\cal O}(\varepsilon^2\delta).$$
It remains to consider the viscous terms on the right hand side of (\ref{eq:NS_2d_new2}).
The first one is similar to the viscous part of the pressure term.
For the second one, using finite differences along the vertical, we write
\begin{eqnarray*}
\int_{z_{\alpha-1/2}}^{z_{\alpha+1/2}} \frac{\partial
}{\partial z}\left(\nu_0\frac{\partial
u}{\partial z}\right) dz & = & \nu_0\left.\frac{\partial u}{\partial
z}\right|_{z_{\alpha+1/2}} - \nu_0\left.\frac{\partial u}{\partial
z}\right|_{z_{\alpha-1/2}},\\
& \approx  & 2\nu\frac{u_{\alpha+1} - u_\alpha }{h_{\alpha+1}+h_\alpha} - 2\nu\frac{u_{\alpha} - u_{\alpha-1}}{h_{\alpha}+h_{\alpha-1}},
\end{eqnarray*}
and relation (\ref{eq:sv_mc3}) follows.
Notice that equations (\ref{eq:sv_mc2}) and (\ref{eq:sv_mc4}) 
are concerned with the evolution of the discharge 
in the lowest and uppest layers, respectively.
The difference between equations (\ref{eq:sv_mc2}) and (\ref{eq:sv_mc4})
and the general equation (\ref{eq:sv_mc3}) comes from the particular form
of the viscous effect at the bottom and at the free surface.
 
Finally we drop the ${\cal O}(\varepsilon^2\delta)$ terms and recovering the variables with
dimension, we obtain the system (\ref{eq:sv_mc1})-(\ref{eq:sv_mc4}).
\end{proof}

%%%%%%%%%%%%%%%%%%%%%%%%%%%%%%%%%%%%%%%%%%%%%%%%%%%%%%%%%%%%%%%%%%%%%%%%%%%%%%%%%%%%%%%%%%%%%%%
\section{Properties of the multilayer system}
\label{sec:properties}
%%%%%%%%%%%%%%%%%%%%%%%%%%%%%%%%%%%%%%%%%%%%%%%%%%%%%%%%%%%%%%%%%%%%%%%%%%%%%%%%%%%%%%%%%%%%%%%
In this paragraph we examine some properties of the model depicted in Proposition \ref{prop:sv_mc}. 
We study its hyperbolicity and we exhibit an energy inequality and a kinetic interpretation of the system.

%%%%%%%%%%%%%%%%%%%%%%%%%%%%%%%%%%%%%%%%%%%%%%%%%%%%%%%%%%%%%%%%%%%%%%%%%%%%%%%%%%%%%%%%%%%%%%%
\subsection{Hyperbolicity}

For the simplicity of the discussion we mainly restrict in this subsection 
to the two-layer version of the multilayer model.

Let us first say some words about the multilayer system (\ref{eq:audusse1})-(\ref{eq:audusse3}) 
introduced by Audusse \cite{audusse}. This non-miscible multilayer system was proved to be non-hyperbolic. 
In the general case the system exhibits complex eigenvalues.
In the very simple case $u_1=u_2=u$ the eigenvalues of the hyperbolic part was shown to be 
equal to the classical barotropic eigenvalues of the monolayer shallow water system $u+\sqrt{gH}$, $u-\sqrt{gH}$
plus a baroclinic eigenvalue $u$ that is concerned with the interface waves.
Nevertheless the system is not hyperbolic since $u$ is a double eigenvalue 
associated to a one-dimensional eigenspace. 
This lack of hyperbolicity may lead to the development of instabilities at the interface \cite{audusse,pares}.
In \cite{audusse} a technical trick is proposed to cure the problem.
Here we can prove the well-posedness of the system.

\begin{prpstn}
The two-layer version of the multilayer 
Saint-Venant system (\ref{eq:sv_mc1})-(\ref{eq:sv_mc4})
is strictly hyperbolic when the total water height is strictly positive.
\end{prpstn}
\begin{proof}
The two-layer version of the multilayer system depicted in Proposition \ref{prop:sv_mc}
stands~-- we denote $u=u_{3/2}$ with $u=u_1$ or $u=u_2$ (see \ref{eq:upwind})~--
\begin{eqnarray*}
&&\frac{\partial H}{\partial t} + l \frac{\partial H u_1}{\partial x} 
+ (1-l) \frac{\partial H u_2}{\partial x}  = 0,\\
&&\frac{\partial H u_1}{\partial t} + \frac{\partial H u_1^2}{\partial x} + 
\frac{g}{2} \frac{\partial H^2}{\partial x} = 
-gH \frac{\partial z_b}{\partial x} + u_{}\left(\frac{\partial H}{\partial t} + \frac{\partial H u_1}{\partial x}\right) 
- H\frac{\partial p^a}{\partial x}
+ \frac{2\nu}{lH}\left( u_2 - u_1 \right) - \tilde{\kappa}(\bar{\bf v},H) u_1,\\
&&\frac{\partial H u_2}{\partial t} + \frac{\partial H u_2^2}{\partial x} + \frac{g}{2} \frac{\partial H^2}{\partial x} = 
-gH\frac{\partial z_b}{\partial x} + u_{}\left(\frac{\partial H}{\partial t} + \frac{\partial H u_2}{\partial x}\right)
- \frac{2\nu}{(1-l)H}\left(u_2 - u_1\right)- H\frac{\partial p^a}{\partial x}.
\end{eqnarray*}
The previous formulation can be written under the quasi-linear form
$$M(X)\frac{\partial X}{\partial t} + A(X)\frac{\partial X}{\partial x}=S(X),$$
with
\begin{eqnarray*}
&&X=\left(\begin{array}{c} H\\ q_1\\ q_2 \end{array}\right),\quad
M(X)=\left(\begin{array}{ccc} 1 & 0 & 0\\ -u_{} & 1 & 0 \\ -u_{} & 0 & 1
\end{array}\right),\quad
A(X)=\left(\begin{array}{ccc} 0 & l & (1-l) \\ gH - u_1^2 & 2u_1 - u_{} & 0 \\ gH - u_2^2 & 0 & 2u_2 - u_{} \end{array}\right),\\
&& S(X)=\left(\begin{array}{c} 0\\ -gH\frac{\partial z_b}{\partial x} + \frac{2\nu}{lH}\left(u_2-u_1) \right)-{\tilde{\kappa}(\bar{\bf v},H)}u_1-H\frac{\partial p^a}{\partial x}\\
\\
-gH\frac{\partial z_b}{\partial x}-\frac{2\nu}{(1-l)H}\left(u_2-u_1\right)-H\frac{\partial p^a}
{\partial x}\end{array}\right),
\end{eqnarray*}
and $q_i=Hu_i$, $i=\{1,2\}$. 

The three eigenvalues of $M^{-1}(X)A(X)$
are the roots of $D(x)=\det(A-xM)=0$ with
$$D(x)=-x\Pi_{i=1}^2 (2u_i-u_{}-x)
-l(2u_2-u_{}-x)(gH-u_1^2+u_{}x) -(1-l)(2u_1-u_{}-x)(gH-u_2^2+u_{}x).$$
Let us fix $H$, $l$, $u_1$ and $u_2$ in $\R$. 
Let us suppose $u_1 < u_2$ with $u_2=u_1 + \gamma^2$. 
We recall that the value of the interface velocity $u$ is taken equal to $u_1$ or $u_2$
following the direction of the exchange of mass between the two layers. 

Let us first suppose that $u=u_1$.
Then we obviously have
$$D(u_1)= -2gHl\gamma^2 < 0, \quad D(-\infty)=+\infty, \quad D(+\infty)=-\infty,$$
and some computations lead to
$$D(\max(u_2=u_1+\gamma^2,u_1+2l\gamma^2)) > 0,$$
since 
$D(u_2)= (1-2l)gH\gamma^2 + l\gamma^6 > 0$ if $l\leq1/2$ and
$D(u_1+2l\gamma^2)= 2l(1-l)(4l-1)\gamma^6> 0$ if $l>1/2$. 
It follows that $D(x)$ has three real and simple eigenvalues.

Let us now suppose that $u=u_2$. Then we have
$$D(u_2)= 2gH(1-l)\gamma^2 > 0, \quad D(-\infty)=+\infty, \quad D(+\infty)=-\infty,$$
and some computations lead to
$$D(\min(u_1=u_2-\gamma^2,u_2-2(1-l)\gamma^2)) < 0,$$
since 
$D(u_1)=(1-2l)gH\gamma^2-(1-l)\gamma^6 < 0$ if $l\geq1/2$ and
$D(u_2-2(1-l)\gamma^2)=2l(1-l)(4l-3)\gamma^6<0$ if $l<1/2$.
Here also $D(x)$ has three real and simple eigenvalues.

The case $u_2<u_1$ is similar and we can conclude 
that the two-layer version of the multilayer system 
depicted in Proposition \ref{prop:sv_mc} is strictly hyperbolic.
Notice that when $u_1=u_2=u$, we find the same baroclinic and
barotropic eigenvalues $u$, $u+\sqrt{gH}$, $u-\sqrt{gH}$ 
as for the nonmiscible multilayer system \cite{audusse},
but they are all simple eigenvalues in this case 
since we consider a system with only three equations. 
\end{proof}

In the case of $N$ layers the matrices $A(X)$ and $M(X)$ can be written
$${\cal A}_{N+1}=\left(\begin{array}{cccccc}
0 & l_1 & \ldots & \ldots & l_N\\
gH-u_1^2 & 2u_1 & \tilde{v}_{1,2} & \ldots & \tilde{v}_{1,N}\\
\vdots & \bar{v}_{2,1} & \ddots & \ddots & \vdots \\
\vdots & \vdots & \ddots & \ddots & \tilde{v}_{N-1,N}\\
gH-u_N^2 & \bar{v}_{N,1} & \ldots & \bar{v}_{N,N-1} & 2u_N
\end{array}\right),$$
with $\bar{v}_{ij}=u_{i-1/2}*l_j/l_i$ and $\tilde{v}_{ij}=u_{i+1/2}*l_j/l_i$, 

$${\cal M}_{N+1}=\left(\begin{array}{ccccc}
1 & 0 & \ldots & \ldots & 0\\
{v}_1 & 1 & \ddots & & \vdots\\
\vdots & 0 & \ddots & \ddots & \vdots\\
\vdots & \vdots & \ddots & 1 & 0\\
{v}_N & 0 & \ldots & 0 & 1
\end{array}\right),$$
with ${v}_{i}=u_{i-1/2}*\sum_{j=1}^{i-1}l_j/l_i+u_{i+1/2}*\sum_{j=i+1}^{N}l_j/l_i$.

We have perfomed various numerical evaluations of the
eigenelements of the matrix ${\cal M}_{N+1}^{-1}{\cal A}_{N+1}$ with numerous choices of the parameters $H$, $u_\alpha$, $u_{\alpha+1/2}$ and $l_\alpha$. All these tests have always shown that the matrix is diagonalizable on $\R$.
In the simple case where all the layers have the same velocity $u$,
the barotropic eigenvalues $u+\sqrt{gH}$ and $u-\sqrt{gH}$ are simple
and the baroclinic eigenvalue $u$ has a multiplicity of $N-1$
but the matrix remains diagonalizable on ${\R}$ and the problem is still well-posed.
%%%%%%%%%%%%%%%%%%%%%%%%%%%%%%%%%%%%%%%%%%%%%%%%%%%%%%%%%%%%%%%%%%%%%%%%%%%%%%%%%%%%%%%%%%%%%%%
\subsection{Energy equality}
The classical Saint-Venant system (\ref{eq:sv1})-(\ref{eq:sv2}) admits an
energy equality \cite{audusse,JSM_DCDS} under the form
\begin{eqnarray}
\frac{\partial E_{sv}}{\partial t} + \frac{\partial}{\partial x}
\left(\bar{u}\bigl(E_{sv}+g\frac{H^2}{2}\bigr) - 4\nu H\bar{u}\frac{\partial \bar{u}}{\partial x}\right) = H\frac{\partial p^a}{\partial t} - 4\nu H\bigl(\frac{\partial \bar{u}}{\partial x}\bigr)^2
 - \frac{\kappa(\bar{{\bf v}},H)}{1+\frac{\kappa(\bar{{\bf v}},H) H}{3 \nu}}\bar{u}^2  + gH\frac{\partial z_b}{\partial t},\label{eq:energy}
\end{eqnarray}
with $E_{sv}=\frac{H\bar{u}^2}{2}+\frac{gH(\eta+z_b)}{2} + Hp^a$. 
Here we have the following result
\begin{prpstn}
For the multilayer Saint-Venant system (\ref{eq:sv_mc1})-(\ref{eq:sv_mc4}), smooth solutions satisfy the energy equality
\begin{eqnarray}
&&\frac{\partial }{\partial t} \left( \sum_{\alpha=1}^N E_{sv,\alpha}^{mc}\right) + \frac{\partial}{\partial x}
\left(\sum_{\alpha=1}^N u_\alpha\left(E_{sv,\alpha}^{mc} + \frac{g}{2}h_\alpha H - 4\nu h_\alpha\frac{\partial u_\alpha}{\partial x}\right) \right) = \nonumber\\
&&\qquad - \kappa(\bar{{\bf v}},H)u_1^2 - \frac{\nu}{h_\alpha}\sum_{\alpha=1}^{N-1} (u_{\alpha+1/2} - u_{\alpha-1/2})^2  - 4\nu \sum_{\alpha=1}^N h_\alpha\bigl(\frac{\partial u_\alpha}{\partial x}\bigr)^2 
+ H\frac{\partial p^a}{\partial t} + gH\frac{\partial z_b}{\partial t},\label{eq:energy_mc}
\end{eqnarray}
with $E_{sv,\alpha}^{mc}=\frac{h_\alpha u_\alpha^2}{2}+\frac{gh_\alpha(\eta+z_b)}{2} + h_\alpha p^a$.
\label{prop:sv_mc_energy}
\end{prpstn}
\begin{proof} 
The proof relies on
classical computations. Starting from (\ref{eq:NS_2d_new2}) with
$u=u^{mc}$, $p=p_{h}$ multiplying it with $u^{mc}$ and integrating over
$[z_{\alpha-1/2},z_{\alpha+1/2}]$ with $1<\alpha<N$ we obtain
\begin{eqnarray}
&&\frac{\partial }{\partial t} E_{sv,\alpha}^{mc} +
\frac{\partial}{\partial x} \left(u_\alpha\left(E_{sv,\alpha}^{mc} + \frac{g}{2}h_\alpha H - 4\nu h\frac{\partial u_\alpha}{\partial x}\right)\right) = \nonumber\\
&&\qquad-\frac{u^2_{\alpha-1/2}}{2}\left( \frac{\partial z_{\alpha-1/2}}{\partial t}
  + u_{\alpha-1/2}\frac{\partial z_{\alpha-1/2}}{\partial x} -
  w_{\alpha-1/2} \right)
+\frac{u^2_{\alpha+1/2}}{2}\left( \frac{\partial z_{\alpha+1/2}}{\partial t}
  + u_{\alpha+1/2}\frac{\partial z_{\alpha+1/2}}{\partial x} - w_{\alpha+1/2} \right)\nonumber\\
&&\qquad-\nu_0 u_{\alpha-1/2}\left. \frac{\partial u^{mc}}{\partial z}\right|_{z_{\alpha-1/2}} + \nu_0 u_{\alpha+1/2}\left. \frac{\partial u^{mc}}{\partial z}\right|_{z_{\alpha+1/2}}\nonumber\\ 
&&\qquad +\nu_0\frac{\partial z_{\alpha-1/2}}{\partial x}u_{\alpha-1/2}\left. \frac{\partial u^{mc}}{\partial x}\right|_{z_{\alpha-1/2}} - \nu_0\frac{\partial z_{\alpha+1/2}}{\partial x}u_{\alpha+1/2}\left. \frac{\partial u^{mc}}{\partial x}\right|_{z_{\alpha+1/2}}\nonumber\\
&&\qquad  - \frac{\nu}{h_\alpha}(u_{\alpha+1/2} - u_{\alpha-1/2})^2  - 4\nu h_\alpha\bigl(\frac{\partial u_\alpha}{\partial x}\bigr)^2 
+ h_\alpha\frac{\partial p^a}{\partial t} + gh_\alpha\frac{\partial z_b}{\partial t},\label{eq:energ_alpha}
\end{eqnarray}
where we have considered for $z \in [z_{\alpha-1/2},z_{\alpha+1/2}]$
$$\frac{\partial u}{\partial z} = \frac{1}{h_\alpha}\left(u_{\alpha+1/2}-u_{\alpha-1/2}\right).$$
An analoguous calculation is valid for $\alpha=1$ and $\alpha=N$. A sum
from $\alpha=1$ to $\alpha=N$ of the equalities (\ref{eq:energ_alpha})
with the boundary conditions (\ref{eq:BC_ap1})-(\ref{eq:BC_ap5}) completes the proof. 
\end{proof}
%%%%%%%%%%%%%%%%%%%%%%%%%%%%%%%%%%%%%%%%%%%%%%%%%%%%%%%%%%%%%%%%%%%%%%%%%%%%%%%%%%%%%%%%%%%%%%%
\subsection{Kinetic interpretation}
\label{sec:kinetic}
For the simulation of a multilayer system several strategies
are possible. Pares {\it et al.} \cite{pares1} consider the full system
and build a specific solver for the two-layer case. Following
the discrete multilayer scheme proposed by Audusse \cite{audusse} we
prefer to exhibit a kinetic formulation of the system obtained in
Proposition \ref{prop:sv_mc}. Indeed kinetic schemes
might be one of the best compromise between accuracy, stability and
efficiency for the resolution of Saint-Venant type equations, see \cite{bristeau,perthame}. 
We refer to the next section for the presentation of the numerical scheme.
Here we focus on the kinetic interpretation of the system.

The kinetic approach consists in using a description of the
microscopic behavior of the system. In this method, fictitious particles are introduced and the equations
are considered at the microscopic scale, where no discontinuities
occur. The process to obtain the kinetic interpretation of the
multilayer model is similar to the one used in
\cite{bristeau} for the monolayer shallow water system. For a given
layer $\alpha$, a distribution function $M_\alpha(x,t,\xi)$ of
fictitious particles with microscopic velocity $\xi$ is introduced
to obtain a linear microscopic kinetic equation equivalent to
the macroscopic model presented in proposition \ref{prop:sv_mc}.

Let us introduce a real function $\chi$ defined on
$\mathbb{R}$, compactly supported and which have the following properties
\begin{equation}
\left\{\begin{array}{l}
\chi(-w) = \chi(w) \geq 0\\
\int_{\mathbb{R}} \chi(w)\ dw = \int_{\mathbb{R}} w^2\chi(w)\ dw = 1.
\end{array}\right.
\label{eq:chi1}
\end{equation}
Now let us construct a density of particles $M_\alpha(x,t,\xi)$
defined by a Gibbs equilibrium: the microscopic density of particles
present at time $t$ in the layers $\alpha$, in the vicinity $\Delta x$ of the abscissa $x$ and with
velocity $\xi$ given by
\begin{equation}
M_\alpha(x,t,\xi) = l_\alpha\frac{H(x,t)}{c} \chi\left(\frac{\xi -
    u_\alpha(x,t)}{c}\right),\qquad \alpha=1,\ldots,N,
\label{eq:Malpha}
\end{equation}
with
$$c^2 = \frac{gH}{2}.$$
Likewise, we define $N_{\alpha+1/2}(x,t,\xi)$ by
\begin{equation}
N_{\alpha+1/2}(x,t,\xi) = G_{\alpha+1/2}(x,t) \ \delta \left(\xi - u_{\alpha+1/2}(x,t)\right),
\qquad \alpha=0,\ldots,N,
\label{eq:N}
\end{equation}
where $\delta$ denotes the Dirac distribution. The quantities $G_{\alpha+1/2}$, $0\leq\alpha\leq N$ represent the
mass exchanges between layers $\alpha$ and $\alpha+1$, they are defined in
(\ref{eq:Q}) and satisfy the conditions (\ref{eq:Qlim}), so $N_{1/2}$ and $N_{N+1/2}$ also satisfy 
\begin{equation}
N_{1/2}(x,t,\xi) = N_{N+1/2}(x,t,\xi) = 0.
\label{eq:Nlim}
\end{equation}

We also introduce the densities $\widetilde{M}_\alpha(x,t,\xi)$ that will be used for the energy equations , they are defined by
\begin{equation*}
\widetilde{M}_\alpha(x,t,\xi) = \frac{gH(x,t)h_\alpha(x,t)}{4c} \chi\left(\frac{\xi - u_\alpha(x,t)}{c}\right).
\end{equation*}
Notice that the introduction of this second family of densities is not needed
when we consider the two dimensional shallow water system. 
Here they take into account some kind of transversal effect at the kinetic level
that is implicitely included into the macroscopic one dimensional shallow water system.
We refer the reader to \cite{bristeau,simeoni} for more details.

With the previous definitions, dropping the viscous,
and friction terms, we write a kinetic representation of the
multilayer Saint-Venant system described in proposition \ref{prop:sv_mc}
and we have the following proposition:
\begin{prpstn}
The functions $(H,u^{mc})$ are strong solutions of the multilayer Saint-Venant system
(\ref{eq:sv_mc1})-(\ref{eq:sv_mc4}) if and only if the set of
equilibria $\{M_\alpha(x,t,\xi)\}_{\alpha=1}^N$
is solution of the kinetic equations
\begin{eqnarray}
\frac{\partial M_\alpha}{\partial t} + \xi \frac{\partial M_\alpha}{\partial x} - \frac{\partial}{\partial x}\left(p^a+gz_b\right)\frac{\partial M_\alpha}{\partial \xi} 
 - N_{\alpha+1/2}(x,t,\xi) + N_{\alpha-1/2}(x,t,\xi)& =& Q_\alpha(x,t,\xi),\label{eq:gibbs}\\
&&\alpha=1,\ldots,N,\nonumber
\end{eqnarray}
 with $\{N_{\alpha+1/2}(x,t,\xi)\}_{\alpha=0}^N$ satisfying
 (\ref{eq:N}),(\ref{eq:Nlim}). The set of equations (\ref{eq:gibbs}) can
also be written under the form
\begin{equation}
\displaystyle
N_{\alpha+1/2}(x,t,\xi) = \sum_{i=1}^\alpha \left(\frac{\partial
    M_i}{\partial t} + \xi \frac{\partial M_i}{\partial
    x} - \frac{\partial}{\partial x}\bigl(p^a + z_b\bigr)\frac{\partial
  M_i}{\partial \xi} - Q_i\right),\qquad \alpha=1,\ldots,N.
\label{eq:gibbsn}
\end{equation}
 The quantities $Q_\alpha(x,t,\xi)$ are ``collision terms''  equals to zero at the
macroscopic level i.e. which satisfy for a.e. values of $(x,t)$
$$\int_{\mathbb{R}} Q_\alpha d\xi =0,\qquad \int_{\mathbb{R}} \xi Q_\alpha d\xi =0.$$
The solution of (\ref{eq:gibbs}),(\ref{eq:gibbsn}) is an entropy solution if additionally
\begin{equation}
\frac{\partial {\widetilde M}_\alpha}{\partial t} + \xi \frac{\partial {\widetilde M}_\alpha}{\partial x}  = {\widetilde Q}_\alpha(x,t,\xi),\quad \alpha=1,\ldots,N,
\end{equation}
with
$$\int_{\mathbb{R}} \left(\frac {\xi^2}{2} Q_\alpha +{\widetilde Q}_\alpha \right) d\xi  \leq 0.$$
\label{prop:kinetic_sv_mc}
\end{prpstn}
\begin{proof} 
As previously we denote $X=(H,q_1,\ldots,q_N)^T$ the vector of unknowns with $q_{\alpha}=l_{\alpha}H u_{\alpha}$.
We introduce $M= (M_1,\ldots,M_N)^T$ and an $(N+1)\times N$ matrix ${\cal K(\xi)}$ 
defined by ${\cal K}_{1,j} = 1$, ${\cal K}_{i+1,j} =
\delta_{i,j}\ \xi$ with $\delta_{i,j}$ the Kronecker symbol. 

Using the definition (\ref{eq:Malpha}) and the properties of the
function $\chi$, we have
\begin{equation}
l_\alpha H(x,t)=\int_{\R}   M_\alpha (x,t,\xi) d\xi,
\label{eq:intM}
\end{equation}
and
\begin{equation}
X(x,t)=\int_{\R} {\cal K(\xi)} \ M (x,t,\xi) d\xi.
\label{eq:intxi}
\end{equation}

The proof is obtained by a simple integration in $\xi$ of the set of equations
(\ref{eq:gibbs}) against the matrix ${\cal K (\xi)}$. First,
an integration in $\xi$ of (\ref{eq:gibbs})
gives the continuity equation (\ref{eq:eq3})  i.e.
\begin{equation*}
\frac{\partial l_\alpha H}{\partial t} + \frac{\partial l_\alpha H u_\alpha}{\partial x} = G_{\alpha+1/2} - G_{\alpha-1/2},
\end{equation*}
and by summation we have (\ref{eq:sv_mc1}). Actually from the
definition (\ref{eq:N}) of $N_{\alpha+1/2}$  we have
$$\int_{\R} N_{\alpha+1/2}(x,t,\xi) d\xi = G_{\alpha+1/2}(x,t),$$
and
$$\int_{\R} \xi N_{\alpha+1/2}(x,t,\xi) d\xi =u_{\alpha+1/2} G_{\alpha+1/2}.$$

Likewise for the
energy balance of the layer $\alpha$ we proceed an integration in $\xi$ of
(\ref{eq:gibbs}) against ${\xi^2}/{2}$. Since we have
\begin{equation}
\int_\R \left(\frac{\xi^2}{2}M_\alpha + \widetilde{M}_\alpha\right) d\xi = \frac{h_\alpha}{2} u_\alpha^2 + \frac{g}{2}h_\alpha H,
\label{eq:energ_term1}
\end{equation}
\begin{equation}
\int_\R \xi \left(\frac{\xi^2}{2}M_\alpha + \widetilde{M}_\alpha\right) d\xi = \frac{h_\alpha}{2} u_\alpha^3 + g h_\alpha H u_\alpha,
\label{eq:energ_term1bis}
\end{equation}
and for the source term
\begin{eqnarray}
\int_\R \frac{\xi^2}{2}\frac{\partial}{\partial x}(p^a+gz_b)\frac{\partial M_\alpha}{\partial \xi} d\xi & = & -\frac{\partial}{\partial x}( p^a + gz_b )h_\alpha u_\alpha\nonumber\\
& = & -\frac{\partial }{\partial x}\bigl((gz_b + p^a)h_\alpha u_\alpha\bigr) + (gz_b+p^a)\frac{\partial h_\alpha u_\alpha
}{\partial x}\nonumber\\
& = & -\frac{\partial }{\partial x}\bigl((gz_b + p^a)h_\alpha u_\alpha\bigr) - \left(gz_b + p^a \right) \frac{\partial h_\alpha}{\partial t}\nonumber\\
& & - \left(gz_b + p^a \right) G_{\alpha+1/2} + \left(gz_b + p^a \right) G_{\alpha-1/2},
\label{eq:energ_term2}
\end{eqnarray}
we obtain the equality
\begin{eqnarray}
&&\frac{\partial}{\partial t}\left(\frac{h_\alpha}{2} u_\alpha^2 + \frac{g}{2}h_\alpha (\eta+z_b) + h_\alpha p^a\right) + \frac{\partial }{\partial x}\left[u_\alpha\left(h_\alpha u_\alpha^2 + \frac{g}{2}h_\alpha H
 + \frac{g}{2}h_\alpha (\eta+z_b) + h_\alpha p^a\right)\right] \nonumber\\
&&\qquad + \frac{u_{\alpha-1/2}^2}{2} G_{\alpha-1/2} - \frac{u_{\alpha+1/2}^2}{2} G_{\alpha+1/2}
 - h_\alpha\frac{\partial p^a}{\partial t} - gh_\alpha\frac{\partial z_b}{\partial t} = \int_\R \left(\frac{\xi^2}{2} Q_\alpha +{\widetilde Q}_\alpha\right)d\xi.
\label{eq:energ_cin}
\end{eqnarray}
The previous relation corresponds to (\ref{eq:energ_alpha}) where
 the viscous and friction terms are neglected. The sum of the 
 equations (\ref{eq:energ_cin}) gives the energy equality for the global
 system and that completes the proof.
\end{proof}
The formulation (\ref{eq:gibbs}) reduces the nonlinear multilayer Saint-Venant system 
to a linear transport system on nonlinear
quantities $\{M_\alpha\}_{\alpha=1}^N$, $\{N_{\alpha+1/2}\}_{\alpha=0}^N$ for which it is easier to find a simple
numerical scheme with good theoretical properties. In the case of a
single layer, for a detailed proof of the kinetic interpretation 
refer to \cite{bristeau} and for the treatment of the source term at this microscopic level
see \cite{simeoni}. Notice that the choice of the function $\chi$ remains quite open at this stage
since several functions satisfy the requested properties.
Following this choice the deduced kinetic scheme will have different properties.
%%%%%%%%%%%%%%%%%%%%%%%%%%%%%%%%%%%%%%%%%%%%%%%%%%%%%%%%%%%%%%%%%%%%%%%%%%%%%%%%%%%%%%%%%%%%%%%
\section{Numerical results}
\label{sec:numeric}
%%%%%%%%%%%%%%%%%%%%%%%%%%%%%%%%%%%%%%%%%%%%%%%%%%%%%%%%%%%%%%%%%%%%%%%%%%%%%%%%%%%%%%%%%%%%%%%
%%%%%%%%%%%%%%%%%%%%%%%%%%%%%%%%%%%%%%%%%%%%%%%%%%%%%%%%%%%%%%%%%%%%%%%%%%%%%%%%%%%%%%%%%%%%%%%
 In the applications discussed here, we assume $p^a=0$ and we neglect the horizontal viscosity.
Then the $N+1$ equations of the multilayer system (\ref{eq:sv_mc1})-(\ref{eq:sv_mc4}) can be written with the general form
\begin{eqnarray}
&&\frac{\partial H}{\partial t} + \sum_{\alpha=1}^N \frac{\partial (l_{\alpha} H u_\alpha)}{\partial x}  = 0,\label{eq:sv_mcg1}\\
\nonumber\\
&&\frac{\partial (l_{\alpha}H u_\alpha)}{\partial t} + \frac{\partial}{\partial x} (l_{\alpha}H u_\alpha^2 
+ \frac{g}{2} l_{\alpha}H^2) = 
-g l_{\alpha}H\frac{\partial z_b}{\partial x}
+u_{\alpha+1/2} G_{\alpha+1/2} -u_{\alpha-1/2} G_{\alpha-1/2}\nonumber\\
&&\qquad + \frac{2\nu_{\alpha}}{l_{\alpha+1}+l_\alpha} \frac{u_{\alpha+1} - u_\alpha }{H} - \frac{2\nu_{\alpha-1}}{l_{\alpha}+l_{\alpha-1}}\frac{u_{\alpha} - u_{\alpha-1}}{H}- \kappa_{\alpha}(u,H) u_{\alpha},
\qquad  \alpha =1,\ldots,N,
\label{eq:sv_mcg}
\end{eqnarray}
with
$$
\kappa_{\alpha}=
\left\{
\begin{array}{l}
\kappa(u,H) \quad {\rm if} \quad \alpha =1\\
0 \quad {\rm if} \quad\alpha \ne 1
\end{array}
\right.\quad
\nu_{\alpha}=
\left\{
\begin{array}{l}
0 \quad {\rm if} \quad \alpha =0\\
\nu \quad {\rm if} \quad\alpha= 1,...,N-1\\
0\quad {\rm if} \quad\alpha=N
\end{array}
\right.$$
  
The previous system is of the form:
\begin{equation}
\frac{\partial X}{\partial t} + \frac{\partial F(X)}{\partial x}=Sb(X)+ Se(X) + Sv (X)
\label{eq:glo}
\end{equation}
 with $F(X)$ the flux of the hyperbolic part, $Sb(X)$ the topography source term, $Se(X)$ the
mass transfer source term and $Sv (X)$ the viscous and friction terms. 

To approximate the solution of the multilayer Saint-Venant system, we use a finite volume framework.
We assume that the computational domain is discretised by $I$ nodes $x_i$.
  We denote $C_i$ the cell of length $\Delta
x_i=x_{i+1/2}-x_{i-1/2}$ with $x_{i+1/2}=(x_i+x_{i+1})/2$. For the time discretization, we denote $t^n = \sum_{k \leq n} \Delta t^k$ where the time steps $\Delta t^k$ will be
precised later though a CFL condition. We denote $X^n_{i}=(H^n_i,q^n_{1,i},\ldots,q^n_{N,i})$ the approximate solution at time $t^n$ on the cell $C_i$ with $q^n_{\alpha,i}=l_{\alpha} H^n_i u^n_{\alpha,i}$.

\subsection{Time discretization}
\label{subsec:timedis}
 For the time discretization, we apply time splitting to the equation (\ref{eq:glo}) and we  write
\begin{eqnarray}
&&\frac{{\tilde X}^{n+1}-X^{n}}{\Delta t^n} + \frac{\partial F(X^n)}{\partial x}=Sb(X^n)+ Se(X^n),
\label{eq:glo1}\\
&&\frac{X^{n+1}-{\tilde X}^{n+1}}{\Delta t^n} -Sv (X^{n},X^{n+1})=0.
\label{eq:glo2}
\end{eqnarray}

Classically we first compute the hyperbolic part (\ref{eq:glo1}) of the multilayer system 
by an explicit scheme. This first computation includes the topographic source term 
in order to preserve relevant equilibria \cite{bristeau1} and also defines the mass transfer terms. 
Concerning the viscous and friction terms (\ref{eq:glo2}) that are dissipative, we prefer a semi-implicit scheme for reasons of stability.

\subsection{Numerical scheme : explicit part}
\label{subsec:scheme}

To perform the explicit step we deduce 
a finite volume kinetic scheme from the previous kinetic
interpretation of the multilayer system. Notice that even if the system is hyperbolic, the eigenvalues
are unknown. Thus any solver requiring the knowledge of the eigenvalues
 while but the kinetic scheme is easily extended \cite{bristeau3}.

Starting from a piecewise constant approximation of the initial
data, the general form of a finite volume method is
\begin{equation}
{\tilde X}^{n+1}_{i}-X^{n}_{i} + \sigma_i^n
\left[F^n_{i+1/2}- F^n_{i-1/2}\right] = \Delta t^n Sb_{i}^{n} + \Delta t^n Se_{i}^{n},
\label{eq:fv}
\end{equation}
where 
$\sigma^n_i=\Delta t^n/\Delta x_i$ is the ratio between space and time
steps and the numerical flux $F^n_{i+1/2}$ is an
approximation of the exact flux estimated at point $x_{i+1/2}$.

The topographic source term $Sb_{i}^{n}$
is not deduced from the kinetic interpretation (see \cite{simeoni}) but computed by hydrostatic reconstruction, see prop.~\ref{prop:hydro_rec}. As in \cite{bristeau, bristeau2} the kinetic interpretation (\ref{eq:gibbs}) is used to precise the
expression of the fluxes $F^n_{i+1/2}$ in (\ref{eq:fv}). 
First, by analogy with (\ref{eq:Malpha}) we define the discrete
densities of particles $M_{\alpha,i}^n$ by	
$$M_{\alpha,i}^n (\xi)= l_{\alpha}\frac{H^n_{i}}{c^n_i} \chi\left(\frac{\xi - u^n_{\alpha,i}}{c^n_i}\right), 
\qquad\mbox{with}\ c^n_i= \sqrt{ \frac{gH^n_{i}}{2}}.$$

Then the equation (\ref{eq:gibbs}) without the atmospheric pressure and topographic terms
is discretised for each $\alpha$ 
by applying a simple upwind scheme for the advection term
\begin{equation}
f_{\alpha,i}^{n+1} (\xi) = M_{\alpha,i}^{n} (\xi) - \xi\sigma^n_i \left(M_{\alpha,i+1/2}^{n} (\xi) - M_{\alpha,i-1/2}^{n} (\xi)\right)+ \Delta t^n \left( N_{\alpha+1/2,i}^{n+1/2} (\xi) - N_{\alpha-1/2,i}^{n+1/2}(\xi)\right),
\label{eq:cindis}
\end{equation}
where
$$ M_{\alpha,i+1/2}^{n} =
\left\{\begin{array}{ll}
M_{\alpha,i}^{n} & \mbox{if } \xi \geq 0\\
M_{\alpha,i+1}^{n} & \mbox{if } \xi < 0
\end{array}\right.$$
and the terms $ N_{\alpha+1/2,i}^{n+1/2}$ will be defined in the following.

We define the vectors $f^{n+1}_i(\xi)=(f^{n+1}_{1,i}(\xi),\ldots,f^{n+1}_{N,i}(\xi))^T$, $M^n_i(\xi)=(M^n_{1,i}(\xi),\ldots,M^n_{N,i}(\xi))^T$. Each new density function $f_{\alpha,i}^{n+1}$ is not an equilibrium but thanks to the
property of the right hand side of $(\ref{eq:gibbs})$, by analogy with (\ref {eq:intM}),(\ref{eq:intxi}) we can recover
the macroscopic quantities at time $t^{n+1}$. 
We write
\begin{equation}
l_\alpha H^{n+1}_{i}=\int_{\R} f^{n+1}_i (\xi) d\xi,
\label{eq:Hdis}
\end{equation}
and by a simple integration in
$\xi$ of (\ref{eq:cindis}) against ${\cal K(\xi)}$, we can precise the macroscopic formula (\ref{eq:fv}) (without the topographic term) 
\begin{equation}
{\tilde X}^{n+1}_{i}=\int_{\R} {\cal K(\xi)}\  f^{n+1}_i (\xi) d\xi.
\label{eq:xn1}
\end{equation}
If we denote
$$F^n_{i+1/2} = F(X^n_i,X^n_{i+1}) = F^+(X^n_i) + F^-(X^n_{i+1}),$$
we define
\begin{equation}
F^-(X^n_i) = \int_{\xi \in {\R}^-}  \xi {\cal K( \xi)} \ M^n_i(\xi)\ d \xi,\qquad
F^+(X^n_i) = \int_{\xi \in {\R}^+}  \xi {\cal K(\xi) } \ M^n_i(\xi)\ d \xi.
\end{equation}
 More precisely the expression of
$F^+(X_i)$ can be written
\begin{equation}
F^+(X_i) =
\left(\begin{array}{c}
F^+_H(X_i)\\
F^+_{q_1}(X_i)\\
\vdots\\
F^+_{q_N}(X_i)
\end{array}\right),
\label{eq:flux}
\end{equation}
with
\begin{eqnarray*}
&&F^+_H(X_i) = \sum_{\alpha=1}^N F^+_{h_\alpha}(X_i) =\sum_{\alpha=1}^N l_{\alpha}H \int_{w\geq -\frac{u_{\alpha,i}}{c_i}}
(u_{\alpha,i} + w c_i)\chi(w)\ dw,\\
&&F^+_{q_\alpha}(X_i) = l_{\alpha} H \int_{w\geq -\frac{u_{\alpha,i}}{c_i}}
(u_{\alpha,i} + w c_i)^2 \chi(w)\ dw.
\end{eqnarray*}
 We denote also
\begin{equation}
{\cal F}_{h_\alpha,i}=F_{h_\alpha,i+1/2}-F_{h_\alpha,i-1/2}=F^+_{h_\alpha}(X_i) + F^-_{h_\alpha}(X_{i+1})
-\left(F^+_{h_\alpha}(X_{i-1}) + F^-_{h_\alpha}(X_{i})\right).
\label{eq:Fhdis}
\end{equation}

This kinetic method is interesting because it gives a very simple and
natural way to propose a numerical flux through the kinetic
interpretation. If we can perform analytically the integration in
(\ref{eq:flux}), i.e. if the probability function 
$\chi$ defined in (\ref{eq:chi1}) is chosen to be simple enough, 
it is also numerically powerfull because the kinetic level
disappears and the scheme is written directly as a macroscopic scheme
for which only very simple computations are needed. In this paper we
have used
$$\chi(w)=\sqrt{\frac{2}{3}} 1_{|w| \leq \sqrt{\frac{3}{2}}}(w).$$

Let us now precise the terms $ N_{\alpha+1/2,i}^{n+1/2}$ and so the exchange terms $Se_i^n$ defined by
\begin{equation}
Se_i^n=\int_{\R} {\cal K(\xi)}\  \left( N_{\alpha+1/2,i}^{n+1/2}(\xi) - N_{\alpha-1/2,i}^{n+1/2}(\xi)\right) d\xi.
\label{eq:se}
\end{equation}
From the conditions (\ref{eq:Nlim}) we prescribe
\begin{equation}
N_{1/2,i}^{n+1/2}(\xi) = N_{N+1/2,i}^{n+1/2} (\xi)= 0.
\label{eq:Nlimdis}
\end{equation}
So we recover $Se_{H,i}^n=0$ and the equation (\ref{eq:xn1}) defines
$H_i^{n+1}$. By summation of (\ref{eq:cindis}) we have
\begin{equation}
 \Delta t^n N_{\alpha+1/2,i}^{n+1/2} (\xi)
=\sum_{j=1}^\alpha \left( f_{j,i}^{n+1} (\xi) - M_{j,i}^{n} (\xi) 
+\xi\sigma^n_i \left(M_{j,i+1/2}^{n} (\xi) - M_{j,i-1/2}^{n} (\xi)\right)\right),
\qquad \alpha=1,\ldots,N-1,
\label{eq:Ndis}
\end{equation}
and we define
\begin{equation}
G_{\alpha+1/2,i}^{n+1/2} =\int_{\R} N_{\alpha+1/2,i}^{n+1/2}(\xi) d\xi,
\qquad \alpha=0,\ldots,N,
\label{eq:Gdis1}
\end{equation}
so we can write
\begin{equation}
\Delta t^n G_{\alpha+1/2,i}^{n+1/2}=\sum_{j=1}^\alpha\left[ l_j (H_i^{n+1}-H_i^{n})
+\sigma^n_i  (  F_{hj,i+1/2}^n- F_{hj,i-1/2}^n)\right], \qquad \alpha=1,\ldots,N. 
\label{eq:Gdis2}
\end{equation}

Then using the discrete mass conservation equation giving $H_i^{n+1}$,
 the terms $G_{\alpha+1/2,i}^{n+1/2}$ can be written under an explicit form (see
(\ref{eq:Qdx})) i.e. depending only of $X^n_i$
\begin{equation}
\Delta x_i G_{\alpha+1/2,i}^{n+1/2}=\sum_{j=1}^\alpha \left( {\cal F}^n_{h_j,i} 
  - l_j \sum_{p=1}^N {\cal  F}_{h_p,i}^n 
     \right),
\label{eq:Gdis3}
\end{equation}
we have to notice that this definition is compatible with the free surface condition of (\ref{eq:Nlimdis}).

We define
\begin{equation}
N_{\alpha+1/2,i}^{n+1/2}(\xi)=G_{\alpha+1/2,i}^{n+1/2}\ \delta \left(\xi - u_{\alpha+1/2,i}^n\right),
\label{eq:Ndis1}
\end{equation}
with, according to (\ref{eq:upwind})
$$ u_{\alpha+1/2,i}^{n} =
\left\{\begin{array}{ll}
u_{\alpha+1,i}^{n} & \mbox{if }\  G_{\alpha+1/2,i}^{n+1/2} \geq 0,\\
u_{\alpha,i}^{n} & \mbox{if } \ G_{\alpha+1/2,i}^{n+1/2} < 0.
\end{array}\right.$$
Then the exchange term $Se_i^n$ in (\ref{eq:se}) is completely defined. 

We have denoted the
approximations in time of $N_{\alpha+1/2}$ and $G_{\alpha+1/2}$ with an upperscript $n+1/2$
because we have to define $H_i^{n+1}$ at the macroscopic level to obtain the microscopic approximation 
of $N_{\alpha+1/2,i}$ which is used for the computation of the momentum $l_\alpha H_i^{n+1}u_{\alpha,i}^{n+1}$.

The source term $Sb^n_{i}=(Sb^n_{H,i},Sb^n_{1,i},\ldots,Sb^n_{N,i})$ is an approximation 
of the topographic source terms. 
For stability purpose, see \cite{bristeau1} 
we use the following discretization
\begin{equation}
Sb^n_{H,i}=0, \qquad Sb^n_{\alpha,i} = l_{\alpha}
\left(\frac{g}{2}(H^n_{i+1/2-})^2 - \frac{g}{2}(H^n_{i-1/2+})^2\right)
\label{eq:source}
\end{equation}
with
\begin{equation}
\begin{split}
z_{b,i+1/2} & = \max\{z_{b,i},z_{b,i+1}\},\\
H^n_{i+1/2-} & = H^n_{i} + z_{b,i} - z_{b,i+1/2},\\
H^n_{i+1/2+} & = H^n_{i+1} + z_{b,i+1} - z_{b,i+1/2}.
\end{split}
\label{eq:source1}
\end{equation}
And we have the following proposition
\begin{prpstn}
The discretization of the source terms given by
(\ref{eq:source}),(\ref{eq:source1}) preserves the steady states
$$\{u_{\alpha,i}^n = 0\}_{\alpha=1}^N,\quad H_{i}^n + z_{b,i} = Cst \in \R,\qquad 
\forall i,\quad\forall n.$$
given by a  ``lake at rest''.
\label{prop:hydro_rec}
\end{prpstn}
\begin{proof}
For the proof of this proposition, the readers can refer to \cite{bristeau1}.
\end{proof}

The scheme explained in this paragraph allows to calculate ${\tilde X}^{n+1}$ given by (\ref{eq:glo1}) and (\ref{eq:fv}).
%%%%%%%%%%%%%%%%%%%%%%%%%%%%%%%%%%%%%%%%%%%%%%%%%%%%%%%%%%%%%%%%%%%%%%%%%%%%%%%%%%%%%%%%%%%%%%%
\subsection{Numerical scheme : implicit part}
Now we aim to calculate ${X}^{n+1}$ from (\ref{eq:glo2}). Neglecting the horizontal viscosity, 
the vertical viscosity source term can be interpreted as a friction term between
one layer and the two adjacent ones. As usual we
treat this friction term implicitly. This leads to solve a linear
system.

The implicit step does not affect the discrete water height therefore
$$H^{n+1}_{i}={\tilde H}^{n+1}_{i},$$
and the computation of the new velocities $\{u_{\alpha,i}^{n+1}\}_{\alpha=1}^N$ leads to
solve a tridiagonal $N \times N$ linear system that reads
$$T^{n,n+1}_i U_i^{n+1}={\tilde q}^{n+1}_i,$$
with $U_i^{n+1}=(u_{1,i}^{n+1},\ldots,u_{N,i}^{n+1})^T$,
${\tilde q}_i^{n+1}=({\tilde q}_{1,i}^{n+1},\ldots,{\tilde q}_{N,i}^{n+1})^T$ and
\begin{eqnarray*}
T^{n,n+1}_i(1,1)\qquad &=&l_1 H^{n+1}_{i}+\frac{2\Delta t^n}{H^{n+1}_{i}}
\left(\frac{\nu_1}{l_1+l_{2}}\right)+\Delta t^n \kappa(X_i^n,H_i^{n+1}),\\ 
T^{n,n+1}_i(\alpha,\alpha) \qquad &=&l_\alpha H^{n+1}_{i}+\frac{2\Delta t^n}{H^{n+1}_{i}}
\left(\frac{\nu_\alpha}{l_\alpha+l_{\alpha+1}}
+\frac{\nu_{\alpha-1}}{l_\alpha+l_{\alpha-1}}\right), \quad \mbox {for } \alpha \in\{2,\ldots,N\},\\
T^{n,n+1}_i(\alpha,\alpha+1)&=&-\frac{2\Delta t^n}{H^{n+1}_{i}} 
\left(\frac{\nu_\alpha}{l_\alpha+l_{\alpha+1}}\right), \quad \mbox {for } \alpha \in\{1,\ldots,N-1\},\\
T^{n,n+1}_i(\alpha-1,\alpha)&=&-\frac{2\Delta t^n}{H^{n+1}_{i}}
\left(\frac{\nu_{\alpha-1}}{l_\alpha+l_{\alpha-1}}\right), \quad \mbox {for } \alpha \in\{2,\ldots,N\}.
\end{eqnarray*}
For the friction at the bottom, several models can be used among which
are Navier, Chezy and Strickler laws.
%%%%%%%%%%%%%%%%%%%%%%%%%%%%%%%%%%%%%%%%%%%%%%%%%%%%%%%%%%%%%%%%%%%%%%%%%%%%%%%%%%%%%%%%%%%%%%%
\subsection{Stability of the scheme}

We now establish the stability property of the kinetic
scheme. Classically for the Saint-Venant system, a CFL condition
ensures the water height is non negative. This CFL condition means that
the quantity of water leaving a given cell during a time step $\Delta t^n$ is less
than the actual water in the cell.

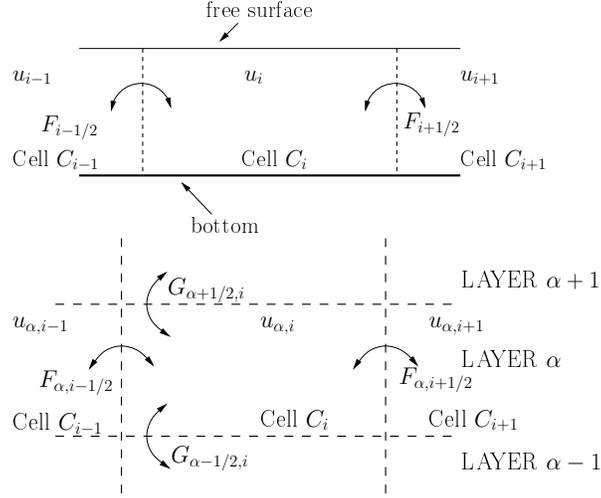
\begin{figure}[htbp]
\begin{center}
 \begin{minipage}{0.8\textwidth}
 \centering
\resizebox{6cm}{!} {\input {Figures/cfl.pstex_t}}
%\resizebox{6cm}{!} {\input {Figures/cfl.pdf_t}}
  \end{minipage}
\vspace{2.cm}
\begin{minipage}{0.8\textwidth}
\centering
\resizebox{6cm}{!}{\input {Figures/cfl_mc.pstex_t}}
%\resizebox{6cm}{!}{\input {Figures/cfl_mc.pdf_t}}
  \end{minipage}
\caption{Interpretation of the CFL condition for the classical
 Saint-Venant system (up) and for the multilayer system (down).}
\label{fig:cfl}
\end{center}
\end{figure}

For the multilayer Saint-Venant system we have the same kind of requirement
concerning the time step $\Delta t^n$. But due to the
vertical discretization, the water can leave the cell $C_i$ of the layer
$\alpha$ either by the
boundaries $x_{i \pm 1/2}$ or by the interfaces $z_{\alpha \pm 1/2}$, see
Fig.~\ref{fig:cfl}. This makes the CFL condition more restrictive
and we have the following proposition

\begin{prpstn}
Assume that the function $\chi$ has a compact support
 of length $2w_M$ then under the CFL condition
\begin{equation}
\Delta t^n \leq \min_{1 \leq \alpha \leq N}\min_{i \in I}
\frac{l_\alpha H_{i}^n\Delta x_i}{l_\alpha H_i^n\left( |u_{\alpha,i}^n| + w_M c_i^n
  \right) + {\Delta x_i}\left( \left[G_{\alpha+1/2,i}^{n+1/2}\right]_- + \left[
      G_{\alpha-1/2,i}^{n+1/2} \right]_+\right)}
\label{eq:cfl}
\end{equation}
the kinetic scheme (\ref{eq:fv}), (\ref{eq:source}) and
(\ref{eq:flux}) keeps the water height
 positive i.e. $H_i^n \geq 0$ if it is true initially. Notice that this
 condition does not depend on $\frac{\partial z_b}{\partial x}$.
\label{prop:positive}
\end{prpstn}
\begin{proof}
The proof has been adapted from those given in \cite{bristeau4,simeoni}. To prove the stability property of the scheme, we come back to the kinetic interpretation and we
 proceed by induction. We assume that $H_{i}^n \geq 0$, $\forall i$ and
 we prove that $H_{i}^{n+1} \geq 0$, $\forall i $.

From the definition of the functions $M_\alpha$ in (\ref{eq:Malpha})
and the positivity of the function $\chi$, we deduce
$$ M_{\alpha,i}^n \geq 0, \qquad \forall i, \quad \mbox {for } \alpha =1,\ldots,N. $$
We now introduce the quantities
$$[\xi]_+ = \max(0,\xi), \quad [\xi]_- = \max(0,-\xi),$$
and so we can write the upwind microscopic scheme
 (\ref{eq:cindis}) 
\begin{eqnarray}
f_{\alpha,i}^{n+1}& = & \left(1 - \sigma_i^n |\xi|\right)M_{\alpha,i}^{n} + \sigma_i^n [\xi]_+ M_{\alpha,i-1}^{n} + \sigma_{i}^n [\xi]_- M_{\alpha,i+1}^{n} \nonumber\\
 & + &\Delta t^n \left(\left(\left[N_{\alpha+1/2,i}^{n+1/2}\right]_+ -\left[N_{\alpha+1/2,i}^{n+1/2}\right]_-\right)
-\left(\left[N_{\alpha-1/2,i}^{n+1/2}\right]_+ -\left[N_{\alpha-1/2,i}^{n+1/2}\right]_-\right)\right).\label{eq:micro}
\end{eqnarray}
The quantity
$$\sigma_i^n |\xi| M_{\alpha,j}^{n} + \Delta t^n \left( \left[N_{\alpha+1/2,i}^{n+1/2}\right]_- + \left[N_{\alpha-1/2,i}^{n+1/2}\right]_+ \right),$$
represents, at the microscopic level, the water leaving the cell $C_i$ of the layer $\alpha$ during
 $\Delta t^n$. A sufficient condition to obtain the stability property, i.e.
\begin{equation}
l_\alpha H_{i}^{n+1} = \int_\R f_{\alpha,i}^{n+1} d\xi \geq 0,
\qquad \forall i, \quad \mbox {for } \alpha =1,\ldots,N,
\label{eq:cond00}
\end{equation}
is then
\begin{equation}
\int_\R \left(\sigma_i^n |\xi| M_{\alpha,i}^{n} + \Delta t^n \left(
    \left[N_{\alpha+1/2,i}^{n+1/2}\right]_- + \left[
      N_{\alpha-1/2,i}^{n+1/2} \right]_+ \right)\right) d\xi \leq
\int_\R M_{\alpha,i}^n d\xi,
\label{eq:cond0}
\end{equation}
and this requirement is satisfied when
$$  \sigma^n_i\left(|u_{\alpha,i}^n| + w_M c_i^n \right) l_\alpha H_{i
}^n + \Delta t^n
\left( \left[G_{\alpha+1/2,i}^{n+1/2}\right]_- + \left[
      G_{\alpha-1/2,i}^{n+1/2} \right]_+\right) \leq l_\alpha H_{i}^{n}.$$
We recall that we have obtained in (\ref{eq:Gdis3}) an explicit form of
$G_{\alpha+1/2,i}^{n+1/2}$. If $ \Delta t^n $ satisfies (\ref{eq:cfl}),
then the condition (\ref{eq:cond00}) is satisfied
and that completes the proof.
\end{proof}

%%%%%%%%%%%%%%%%%%%%%%%%%%%%%%%%%%%%%%%%%%%%%%%%%%%%%%%%%%%%%%%%%%%%%%%%%%%%%%%%%%%%%%%%%%%%%%%
\subsection{Second order scheme}

The second-order accuracy in time is usually recovered by the Heun
method \cite{heun} that is a slight modification of the second order
Runge-Kutta method. The advantage of the Heun scheme is that it
preserves the invariant domains without any additional limitation on the CFL.

We also apply a formally second order scheme in space by a limited reconstruction of the
variables. An advantage of the new multilayer approach with only one continuity equation is
that the water height can be reconstructed while preserving the mass conservation 
without difficulty. 
%%%%%%%%%%%%%%%%%%%%%%%%%%%%%%%%%%%%%%%%%%%%%%%%%%%%%%%%%%%%%%%%%%%%%%%%%%%%%%%%%%%%%%%%%%%%%%%
\subsection{Numerical simulations}
%%%%%%%%%%%%%%%%%%%%%%%%%%%%%%%%%%%%%%%%%%%%%%%%%%%%%%%%%%%%%%%%%%%%%%%%%%%%%%%%%%%%%%%%%%%%%%%
\subsubsection{Transcritical flow over a bump}

We first consider an academic test case that is very commonly used 
for the validation of classical one-layer shallow water solvers.
Here we add some friction at the bottom in order to compare
solutions of one-layer and multilayer shallow water systems
with the solution of hydrostatic incompressible Navier-Stokes equations.
We impose an inflow (left boundary) of $1.0\ m^2.s^{-1}$ 
and the water height at the exit (right boundary) is prescribed to be equal to $0.6\ m$. 
The Strickler friction coefficient at the bottom is $30\ m^{1/3}.s^{-1}$ 
and the kinematic viscosity is $0.01\ m^2.s^{-1}$.
The data are chosen such that the flow is supposed to reach a stationary regime 
that presents some transitions between sub- and supercritical parts and an hydraulic jump.
Notice that an analytical solution exists for this test in the case of a single layer \cite{bristeau,simeoni}. 

\begin{center}
\begin{figure}[htbp]
\includegraphics[height=08cm]{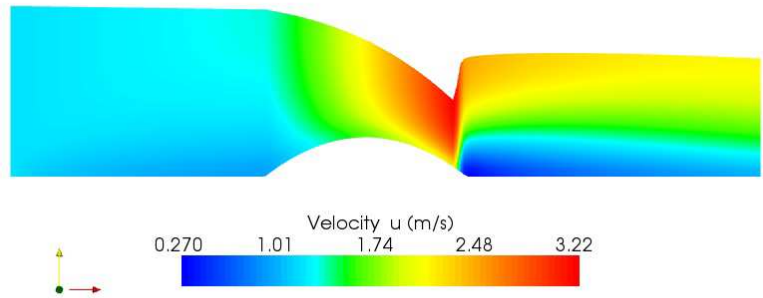}
\caption{Horizontal velocities $\{u_\alpha(x,t_f)\}_{\alpha=1}^N$ with $N=15$ layers.}
\label{fig:u15}
\end{figure}
\end{center}

\begin{center}
\begin{figure}[htbp]
\includegraphics[height=9cm,width=11.5cm]{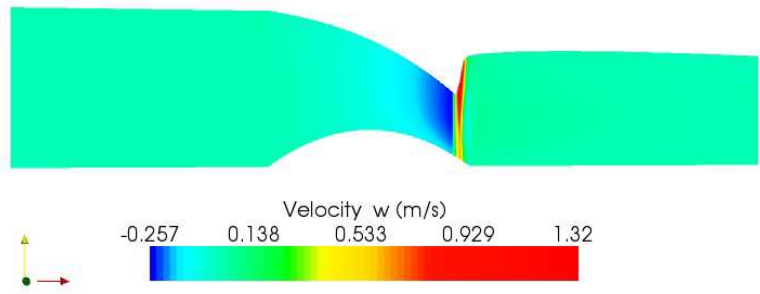}
\caption{Vertictal velocity $\{w_\alpha(x,t_f)\}_{\alpha=1}^N$ with $N=15$ layers.}
\label{fig:w15}
\end{figure}
\end{center}

\begin{center}
\begin{figure}[htbp]
\includegraphics[height=8cm]{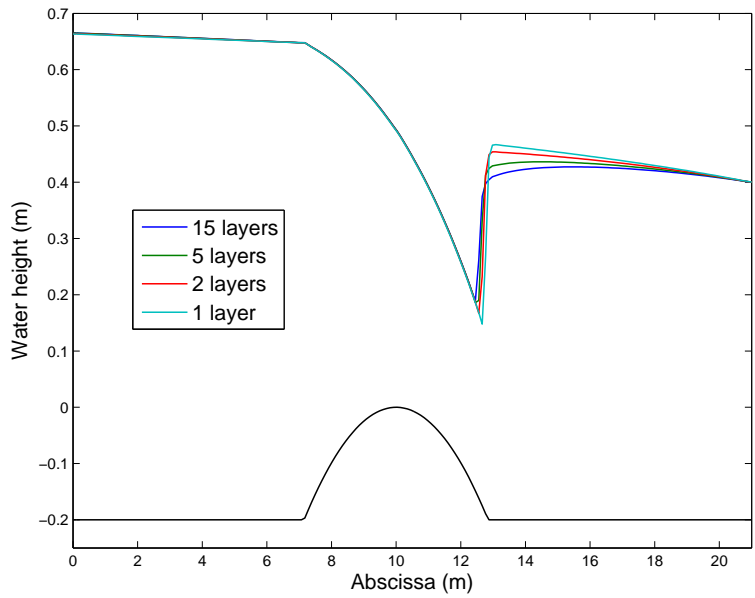}
\caption{Shape of the free surface for simulations carried out with
  different number of layers.}
\label{fig:compar_mc}
\end{figure}
\end{center}

The simulation results are depicted in Fig.~\ref{fig:u15}, \ref{fig:w15} and \ref{fig:compar_mc}. 
The presented results correspond to a time instant $t_f$ where the permanent regime is achieved.
Notice that we present some results related to the vertical velocity in Fig.~\ref{fig:w15}.
Since we consider a shallow water type system we do not need this vertical velocity for the computation.
But it is possible to recover it for postprocessing purpose :
departing from the computed horizontal velocity 
we use the divergence free condition (\ref{eq:NS_2d1}) 
and the non penetration condition at the bottom (\ref{eq:bottom1})
to evaluate an approximation of the vertical velocity.
Notice also that the actual computations are purely one dimensional.
Hence Fig.~\ref{fig:u15} and \ref{fig:w15} present velocity results on a postprocessing mesh that 
is constructed departing from the 1d mesh by the use of the computed layer water heights.

The results depicted in Fig.~\ref{fig:u15} and \ref{fig:w15} 
are consistent with computations performed using the hydrostatic Navier-Stokes equations
\cite{bristeau2} 
and also using the former multilayer Saint-Venant system \cite{audusse}. 
The results depicted in Fig.~\ref{fig:compar_mc} exhibit that the presented solver is quite robust since it
is able to compute transcritical solutions and shock waves even when a large number of layers are considered.
Notice also that the hydraulic jump appears to be overestimated by the one-layer computation
when compared with other results - see Fig.~\ref{fig:compar_mc}. 

%%%%%%%%%%%%%%%%%%%%%%%%%%%%%%%%%%%%%%%%%%%%%%%%%%%%%%%%%%%%%%%%%%%%%%%%%%%%%%%%%%%%%%%%%%%%%%%
\subsubsection{Wind effects}

We claim in the introduction that the great interest of the new multilayer formulation 
that we proposed here is to allow mass exchanges between layers.
This effect is exhibited in the numerical test that we present now.
We consider a lake with a non trivial bottom and vertical shores.
We impose a constant wind stress (from left to right) at the free surface.
The flow is then supposed to reach a stationnary state
that includes some water recirculations in the lake.
Notice that this kind of stationnary flows is clearly
impossible to compute with the classical one-layer shallow water system
since the velocity is imposed to be constant along the vertical. 
They are also out of the domain of application of the former multilayer shallow water system
that was introduced by Audusse \cite{bristeau2} since they clearly involve
large mass transfers (at least near the shores) between the layers. 

\begin{center}
\begin{figure}[htbp]
\includegraphics[height=04cm]{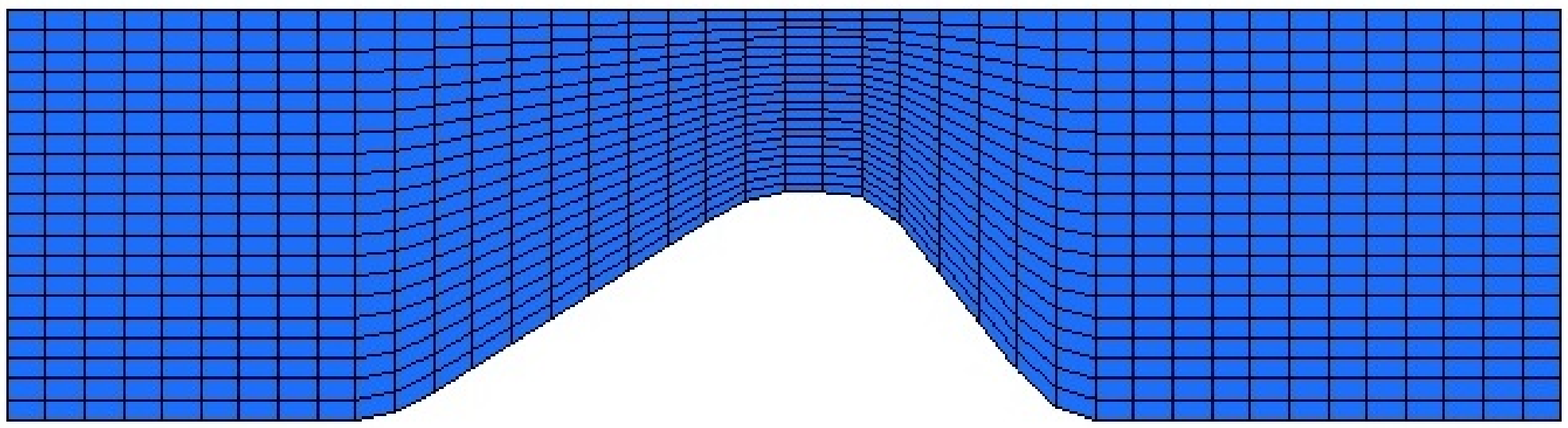}
\caption{The geometrical model with the horizontal mesh and the
  vertical discretization by layers.}
\label{fig:step_wind_mesh}
\end{figure}
\end{center}

\begin{center}
\begin{figure}[htbp]
\includegraphics[height=06cm]{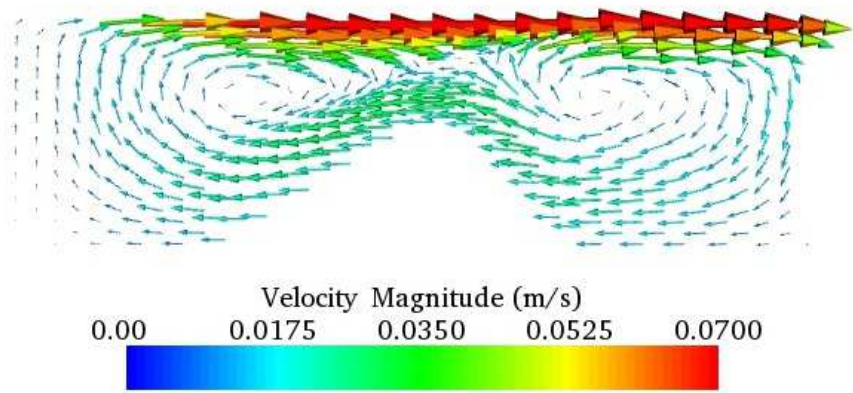}
\caption{A wind blow from the left part of the domain to the right part. The arrows represent the velocity field in the lake.}
\label{fig:step_wind}
\end{figure}
\end{center}

As for the previous case 
we use a reconstruction strategy in order to estimate a vertical velocity field
and we present the results on a postprocessing 2D mesh 
that is presented in Fig.~\ref{fig:step_wind_mesh}.
In Fig.~\ref{fig:step_wind} we present the two dimensional velocity vectors on this 2D mesh.
The results exhibit a global recirculation that is combined with two local recirculations 
that are induced by the topography of the lake. 
The qualitative aspect of the solution is consistent with the previsions.

%%%%%%%%%%%%%%%%%%%%%%%%%%%%%%%%%%%%%%%%%%%%%%%%%%%%%%%%%%%%%%%%%%%%%%%%%%%%%%%%%%%%%%%%%%%%%%%
\section{Conclusion}

In this paper, the authors have described an exchanging mass multilayer
Saint-Venant system. The derivation of the model, the study of its main
properties and a numerical scheme for its discretization are given. Some
simulations are also presented. Notice that the model and the results
presented here in 2D $(x,z)$ are also available in 3D $(x,y,z)$.

Because of its accuracy and simplicity, the kinetic scheme seems well adapted for the simulations of such a
model. Moreover since the eigenvalues of the hyperbolic system are not
explictly known, a lot of finite volume schemes fails in this situation.
 
As depicted in Fig.~\ref{fig:step_wind_mesh}, the vertical
discretization proposed for water height leads to a regular mesh. A
strategy of ``mesh refinement'' based on a inhomogenous number of layers
have to be added.

The presented system can be enriched in several ways. First, the
hydrostatic assumption concerning the pressure terms can be relaxed
leading to the models presented in \cite{JSM_DCDS}. Then we can also
consider a passive pollutant in the flow. This implies to add a
conservation equation for the pollutant concentration. Finally, we can
consider the density of the fluid varies with the concentration of pollutant.
These three improvements have been added to their model by the
authors and will be presented in forthcoming papers.
%%%%%%%%%%%%%%%%%%%%%%%%%%%%%%%%%%%%%%%%%%%%%%%%%%%%%%%%%%%%%%%%%%%%%%%%%%%%%%%%%%%%%%%%%%%%%%%

%\newpage
\bibliographystyle{amsplain}
\bibliography{boussinesq}

\end{document}

%% file: Figures/models_mc.pstex_t
\begin{picture}(0,0)%
\includegraphics{models_mc.pstex}%
\end{picture}%
\setlength{\unitlength}{4144sp}%
\begingroup\makeatletter\ifx\SetFigFont\undefined%
\gdef\SetFigFont#1#2#3#4#5{%
  \reset@font\fontsize{#1}{#2pt}%
  \fontfamily{#3}\fontseries{#4}\fontshape{#5}%
  \selectfont}%
\fi\endgroup%
\begin{picture}(10967,4873)(976,-5777)
\put(2566,-3976){\makebox(0,0)[lb]{\smash{{\SetFigFont{17}{20.4}{\rmdefault}{\mddefault}{\updefault}{\color[rgb]{0,0,0}Hydrostatic models}%
}}}}
\put(4276,-1726){\makebox(0,0)[lb]{\smash{{\SetFigFont{17}{20.4}{\rmdefault}{\mddefault}{\updefault}{\color[rgb]{0,0,0}for incompressible free surface flows}%
}}}}
\put(4816,-2626){\makebox(0,0)[lb]{\smash{{\SetFigFont{17}{20.4}{\rmdefault}{\mddefault}{\updefault}{\color[rgb]{0,0,0}Shallow water assumption}%
}}}}
\put(4951,-1321){\makebox(0,0)[lb]{\smash{{\SetFigFont{17}{20.4}{\rmdefault}{\mddefault}{\updefault}{\color[rgb]{0,0,0}Navier-Stokes equations}%
}}}}
\put(7381,-3976){\makebox(0,0)[lb]{\smash{{\SetFigFont{17}{20.4}{\rmdefault}{\mddefault}{\updefault}{\color[rgb]{0,0,0}Non hydrostatic models}%
}}}}
\put(991,-5011){\makebox(0,0)[lb]{\smash{{\SetFigFont{17}{20.4}{\rmdefault}{\mddefault}{\updefault}{\color[rgb]{0,0,0}Saint-Venant system}%
}}}}
\put(1171,-5461){\makebox(0,0)[lb]{\smash{{\SetFigFont{17}{20.4}{\rmdefault}{\mddefault}{\updefault}{\color[rgb]{0,0,0}Multilayer Saint-Venant system}%
}}}}
\put(5986,-4921){\makebox(0,0)[lb]{\smash{{\SetFigFont{17}{20.4}{\rmdefault}{\mddefault}{\updefault}{\color[rgb]{0,0,0}Boussinesq system}%
}}}}
\put(7021,-5686){\makebox(0,0)[lb]{\smash{{\SetFigFont{17}{20.4}{\rmdefault}{\mddefault}{\updefault}{\color[rgb]{0,0,0}Multilayer extended Saint-Venant system}%
}}}}
\put(6211,-5326){\makebox(0,0)[lb]{\smash{{\SetFigFont{17}{20.4}{\rmdefault}{\mddefault}{\updefault}{\color[rgb]{0,0,0}Extended Saint-Venant system}%
}}}}
\end{picture}%

%% file: Figures/notations_mc.pstex_t
\begin{picture}(0,0)%
\includegraphics{notations_mc.pstex}%
\end{picture}%
\setlength{\unitlength}{4144sp}%
\begingroup\makeatletter\ifx\SetFigFont\undefined%
\gdef\SetFigFont#1#2#3#4#5{%
  \reset@font\fontsize{#1}{#2pt}%
  \fontfamily{#3}\fontseries{#4}\fontshape{#5}%
  \selectfont}%
\fi\endgroup%
\begin{picture}(11118,5238)(-284,-5608)
\put(10531,-1726){\makebox(0,0)[lb]{\smash{{\SetFigFont{20}{24.0}{\rmdefault}{\mddefault}{\updefault}{\color[rgb]{0,0,0}$x$}%
}}}}
\put(1036,-601){\makebox(0,0)[lb]{\smash{{\SetFigFont{20}{24.0}{\rmdefault}{\mddefault}{\updefault}{\color[rgb]{0,0,0}$z$}%
}}}}
\put(3691,-2311){\makebox(0,0)[lb]{\smash{{\SetFigFont{20}{24.0}{\rmdefault}{\mddefault}{\updefault}{\color[rgb]{0,0,0}$\eta(x,t)$}%
}}}}
\put(5716,-1276){\makebox(0,0)[lb]{\smash{{\SetFigFont{20}{24.0}{\rmdefault}{\mddefault}{\updefault}{\color[rgb]{0,0,0}Free surface}%
}}}}
\put(5986,-5371){\makebox(0,0)[lb]{\smash{{\SetFigFont{20}{24.0}{\rmdefault}{\mddefault}{\updefault}{\color[rgb]{0,0,0}Bottom}%
}}}}
\put(8191,-2176){\makebox(0,0)[lb]{\smash{{\SetFigFont{14}{16.8}{\rmdefault}{\mddefault}{\updefault}{\color[rgb]{0,0,0}$h_4(x,t)$}%
}}}}
\put(2926,-4741){\makebox(0,0)[lb]{\smash{{\SetFigFont{20}{24.0}{\rmdefault}{\mddefault}{\updefault}{\color[rgb]{0,0,0}$z_b(x,t)$}%
}}}}
\put(8191,-2896){\makebox(0,0)[lb]{\smash{{\SetFigFont{14}{16.8}{\rmdefault}{\mddefault}{\updefault}{\color[rgb]{0,0,0}$h_3(x,t)$}%
}}}}
\put(8191,-3706){\makebox(0,0)[lb]{\smash{{\SetFigFont{14}{16.8}{\rmdefault}{\mddefault}{\updefault}{\color[rgb]{0,0,0}$h_2(x,t)$}%
}}}}
\put(8191,-4426){\makebox(0,0)[lb]{\smash{{\SetFigFont{14}{16.8}{\rmdefault}{\mddefault}{\updefault}{\color[rgb]{0,0,0}$h_1(x,t)$}%
}}}}
\put(4681,-3436){\makebox(0,0)[lb]{\smash{{\SetFigFont{20}{24.0}{\rmdefault}{\mddefault}{\updefault}{\color[rgb]{0,0,0}$H(x,t)$}%
}}}}
\put(6301,-2131){\makebox(0,0)[lb]{\smash{{\SetFigFont{14}{16.8}{\rmdefault}{\mddefault}{\updefault}{\color[rgb]{0,0,0}$u_4(x,t)$}%
}}}}
\put(6301,-2896){\makebox(0,0)[lb]{\smash{{\SetFigFont{14}{16.8}{\rmdefault}{\mddefault}{\updefault}{\color[rgb]{0,0,0}$u_3(x,t)$}%
}}}}
\put(6301,-3571){\makebox(0,0)[lb]{\smash{{\SetFigFont{14}{16.8}{\rmdefault}{\mddefault}{\updefault}{\color[rgb]{0,0,0}$u_2(x,t)$}%
}}}}
\put(6301,-4291){\makebox(0,0)[lb]{\smash{{\SetFigFont{14}{16.8}{\rmdefault}{\mddefault}{\updefault}{\color[rgb]{0,0,0}$u_1(x,t)$}%
}}}}
\put(1081,-2131){\makebox(0,0)[lb]{\smash{{\SetFigFont{20}{24.0}{\rmdefault}{\mddefault}{\updefault}{\color[rgb]{0,0,0}$0$}%
}}}}
\put(181,-2626){\makebox(0,0)[lb]{\smash{{\SetFigFont{14}{16.8}{\rmdefault}{\mddefault}{\updefault}{\color[rgb]{0,0,0}$z_{3+1/2}(x,t)$}%
}}}}
\put(136,-3346){\makebox(0,0)[lb]{\smash{{\SetFigFont{14}{16.8}{\rmdefault}{\mddefault}{\updefault}{\color[rgb]{0,0,0}$z_{2+1/2}(x,t)$}%
}}}}
\put(136,-4066){\makebox(0,0)[lb]{\smash{{\SetFigFont{14}{16.8}{\rmdefault}{\mddefault}{\updefault}{\color[rgb]{0,0,0}$z_{1+1/2}(x,t)$}%
}}}}
\put(-134,-4831){\makebox(0,0)[lb]{\smash{{\SetFigFont{14}{16.8}{\rmdefault}{\mddefault}{\updefault}{\color[rgb]{0,0,0}$z_{1/2}=z_b(x,t)$}%
}}}}
\put(-269,-1681){\makebox(0,0)[lb]{\smash{{\SetFigFont{14}{16.8}{\rmdefault}{\mddefault}{\updefault}{\color[rgb]{0,0,0}$z_{4+1/2}=\eta(x,t)$}%
}}}}
\end{picture}%

%% file: Figures/cfl.pstex_t
\begin{picture}(0,0)%
\includegraphics{cfl.pstex}%
\end{picture}%
\setlength{\unitlength}{4144sp}%
\begingroup\makeatletter\ifx\SetFigFont\undefined%
\gdef\SetFigFont#1#2#3#4#5{%
  \reset@font\fontsize{#1}{#2pt}%
  \fontfamily{#3}\fontseries{#4}\fontshape{#5}%
  \selectfont}%
\fi\endgroup%
\begin{picture}(6393,3315)(1741,-6736)
\put(2161,-5281){\makebox(0,0)[lb]{\smash{{\SetFigFont{20}{24.0}{\rmdefault}{\mddefault}{\updefault}{\color[rgb]{0,0,0}$F_{i-1/2}$}%
}}}}
\put(4276,-6721){\makebox(0,0)[lb]{\smash{{\SetFigFont{20}{24.0}{\rmdefault}{\mddefault}{\updefault}{\color[rgb]{0,0,0}bottom}%
}}}}
\put(4501,-3661){\makebox(0,0)[lb]{\smash{{\SetFigFont{20}{24.0}{\rmdefault}{\mddefault}{\updefault}{\color[rgb]{0,0,0}free surface}%
}}}}
\put(4996,-5776){\makebox(0,0)[lb]{\smash{{\SetFigFont{20}{24.0}{\rmdefault}{\mddefault}{\updefault}{\color[rgb]{0,0,0}Cell $C_i$}%
}}}}
\put(8101,-4561){\makebox(0,0)[lb]{\smash{{\SetFigFont{20}{24.0}{\rmdefault}{\mddefault}{\updefault}{\color[rgb]{0,0,0}$u_{i+1}$}%
}}}}
\put(8101,-5776){\makebox(0,0)[lb]{\smash{{\SetFigFont{20}{24.0}{\rmdefault}{\mddefault}{\updefault}{\color[rgb]{0,0,0}Cell $C_{i+1}$}%
}}}}
\put(5041,-4561){\makebox(0,0)[lb]{\smash{{\SetFigFont{20}{24.0}{\rmdefault}{\mddefault}{\updefault}{\color[rgb]{0,0,0}$u_i$}%
}}}}
\put(1756,-4561){\makebox(0,0)[lb]{\smash{{\SetFigFont{20}{24.0}{\rmdefault}{\mddefault}{\updefault}{\color[rgb]{0,0,0}$u_{i-1}$}%
}}}}
\put(1756,-5776){\makebox(0,0)[lb]{\smash{{\SetFigFont{20}{24.0}{\rmdefault}{\mddefault}{\updefault}{\color[rgb]{0,0,0}Cell $C_{i-1}$}%
}}}}
\put(7291,-5236){\makebox(0,0)[lb]{\smash{{\SetFigFont{20}{24.0}{\rmdefault}{\mddefault}{\updefault}{\color[rgb]{0,0,0}$F_{i+1/2}$}%
}}}}
\end{picture}%

%% file: Figures/cfl_mc.pstex_t
\begin{picture}(0,0)%
\includegraphics{cfl_mc.pstex}%
\end{picture}%
\setlength{\unitlength}{4144sp}%
\begingroup\makeatletter\ifx\SetFigFont\undefined%
\gdef\SetFigFont#1#2#3#4#5{%
  \reset@font\fontsize{#1}{#2pt}%
  \fontfamily{#3}\fontseries{#4}\fontshape{#5}%
  \selectfont}%
\fi\endgroup%
\begin{picture}(6150,3644)(2101,-6833)
\put(7786,-5776){\makebox(0,0)[lb]{\smash{{\SetFigFont{20}{24.0}{\rmdefault}{\mddefault}{\updefault}{\color[rgb]{0,0,0}Cell $C_{i+1}$}%
}}}}
\put(8236,-4966){\makebox(0,0)[lb]{\smash{{\SetFigFont{20}{24.0}{\rmdefault}{\mddefault}{\updefault}{\color[rgb]{0,0,0}LAYER $\alpha$}%
}}}}
\put(8236,-3886){\makebox(0,0)[lb]{\smash{{\SetFigFont{20}{24.0}{\rmdefault}{\mddefault}{\updefault}{\color[rgb]{0,0,0}LAYER $\alpha+1$}%
}}}}
\put(8236,-6361){\makebox(0,0)[lb]{\smash{{\SetFigFont{20}{24.0}{\rmdefault}{\mddefault}{\updefault}{\color[rgb]{0,0,0}LAYER $\alpha-1$}%
}}}}
\put(4231,-3976){\makebox(0,0)[lb]{\smash{{\SetFigFont{20}{24.0}{\rmdefault}{\mddefault}{\updefault}{\color[rgb]{0,0,0}$G_{\alpha+1/2,i}$}%
}}}}
\put(4276,-6271){\makebox(0,0)[lb]{\smash{{\SetFigFont{20}{24.0}{\rmdefault}{\mddefault}{\updefault}{\color[rgb]{0,0,0}$G_{\alpha-1/2,i}$}%
}}}}
\put(5491,-4426){\makebox(0,0)[lb]{\smash{{\SetFigFont{20}{24.0}{\rmdefault}{\mddefault}{\updefault}{\color[rgb]{0,0,0}$u_{\alpha,i}$}%
}}}}
\put(5491,-5776){\makebox(0,0)[lb]{\smash{{\SetFigFont{20}{24.0}{\rmdefault}{\mddefault}{\updefault}{\color[rgb]{0,0,0}Cell $C_i$}%
}}}}
\put(2476,-5236){\makebox(0,0)[lb]{\smash{{\SetFigFont{20}{24.0}{\rmdefault}{\mddefault}{\updefault}{\color[rgb]{0,0,0}$F_{\alpha,i-1/2}$}%
}}}}
\put(7381,-5191){\makebox(0,0)[lb]{\smash{{\SetFigFont{20}{24.0}{\rmdefault}{\mddefault}{\updefault}{\color[rgb]{0,0,0}$F_{\alpha,i+1/2}$}%
}}}}
\put(2116,-5821){\makebox(0,0)[lb]{\smash{{\SetFigFont{20}{24.0}{\rmdefault}{\mddefault}{\updefault}{\color[rgb]{0,0,0}Cell $C_{i-1}$}%
}}}}
\put(2116,-4426){\makebox(0,0)[lb]{\smash{{\SetFigFont{20}{24.0}{\rmdefault}{\mddefault}{\updefault}{\color[rgb]{0,0,0}$u_{\alpha,i-1}$}%
}}}}
\put(7786,-4426){\makebox(0,0)[lb]{\smash{{\SetFigFont{20}{24.0}{\rmdefault}{\mddefault}{\updefault}{\color[rgb]{0,0,0}$u_{\alpha,i+1}$}%
}}}}
\end{picture}%

%% file: sv+_mc.bbl
\providecommand{\bysame}{\leavevmode\hbox to3em{\hrulefill}\thinspace}
\providecommand{\MR}{\relax\ifhmode\unskip\space\fi MR }
\providecommand{\MRhref}[2]{%
  \href{http://www.ams.org/mathscinet-getitem?mr=#1}{#2}
}
\providecommand{\href}[2]{#2}
\begin{thebibliography}{10}

\bibitem{audusse}
E.~Audusse, \emph{{A multilayer Saint-Venant System~: Derivation and Numerical
  Validation}}, Discrete Contin. Dyn. Syst. Ser. B \textbf{5} (2005), no.~2,
  189--214.

\bibitem{bristeau1}
E.~Audusse, F.~Bouchut, M.O. Bristeau, R.~Klein, and B.~Perthame, \emph{A fast
  and stable well-balanced scheme with hydrostatic reconstruction for {Shallow
  Water} flows}, SIAM J. Sci. Comput. \textbf{25} (2004), no.~6, 2050--2065.

\bibitem{bristeau4}
E.~Audusse and M.O. Bristeau, \emph{Transport of pollutant in shallow water
  flows : A two time steps kinetic method}, M2AN \textbf{37} (2003), no.~2,
  389--416.

\bibitem{bristeau}
\bysame, \emph{A well-balanced positivity preserving second-order scheme for
  shallow water flows on unstructured meshes.}, J. Comput. Phys. \textbf{206}
  (2005), no.~1, 311--333.

\bibitem{bristeau3}
\bysame, \emph{Finite-volume solvers for a multilayer saint-venant system},
  Int. J. Appl. Math. Comput. Sci. \textbf{17} (2007), no.~3, 311--319.

\bibitem{bristeau2}
E.~Audusse, M.O. Bristeau, and Decoene A., \emph{Numerical simulations of 3d
  free surface flows by a multilayer {Saint-Venant} model}, Internat. J. Numer.
  Methods Fluids \textbf{56} (2008), no.~3, 331--350.

\bibitem{saint-venant}
A.J.C. Barr\'e~de Saint-Venant, \emph{{Th\'eorie du mouvement non permanent des
  eaux avec applications aux crues des rivi\`eres et \`a l'introduction des
  mar\'ees dans leur lit}}, C. R. Acad. Sci. Paris \textbf{73} (1871),
  147--154.

\bibitem{heun}
F.~Bouchut, \emph{An introduction to finite volume methods for hyperbolic
  conservation laws.}, ESAIM Proc. \textbf{15} (2004), 107--127.

\bibitem{bouchut1}
F.~Bouchut and T.~Morales~de Luna, \emph{An entropy satisfying scheme for
  two-layer shallow water equations with uncoupled treatment}, M2AN Math.
  Model. Numer. Anal. \textbf{42} (2008), 683--698.

\bibitem{bouchut}
F.~Bouchut and M.~Westdickenberg, \emph{{Gravity driven shallow water models
  for arbitrary topography}}, Comm. in Math. Sci. \textbf{2} (2004), 359--389.

\bibitem{JSM_DCDS}
M.O. Bristeau and J.~Sainte-Marie, \emph{{Derivation of a non-hydrostatic
  shallow water model; Comparison with Saint-Venant and Boussinesq systems}},
  Discrete Contin. Dyn. Syst. Ser. B \textbf{10} (2008), no.~4, 733--759.

\bibitem{pares}
M.J. Castro, J.A. Garc\'{\i}a-Rodr\'{\i}guez, J.M. Gonz\'{a}lez-Vida,
  J.~Mac\'{\i}as, C.~Par\'{e}s, and M.E. V\'{a}zquez-Cend\'{o}n,
  \emph{Numerical simulation of two-layer shallow water flows through channels
  with irregular geometry}, J. Comput. Phys. \textbf{195} (2004), no.~1,
  202--235.

\bibitem{pares1}
M.J. Castro, J.~Mac\'{\i}as, and C.~Par\'{e}s, \emph{A q-scheme for a class of
  systems of coupled conservation laws with source term. application to a
  two-layer {1-D} shallow water system.}, M2AN Math. Model. Numer. Anal.
  \textbf{35} (2001), no.~1, 107--127.

\bibitem{decoene}
A.~Decoene, L.~Bonaventura, E.~Miglio, and F.~Saleri, \emph{Asymptotic
  derivation of the section-averaged shallow water equations for river
  hydraulics}, MOX-Report \textbf{17} (2007).

\bibitem{saleri}
S.~Ferrari and F.~Saleri, \emph{{A new two-dimensional Shallow Water model
  including pressure effects and slow varying bottom topography}}, M2AN Math.
  Model. Numer. Anal. \textbf{38} (2004), no.~2, 211--234.

\bibitem{gerbeau}
J.-F. Gerbeau and B.~Perthame, \emph{{Derivation of Viscous Saint-Venant System
  for Laminar Shallow Water; Numerical Validation}}, Discrete Contin. Dyn.
  Syst. Ser. B \textbf{1} (2001), no.~1, 89--102.

\bibitem{lions}
P.L. Lions, \emph{{Mathematical Topics in Fluid Mechanics. Vol.~1:
  Incompressible models.}}, Oxford University Press, 1996.

\bibitem{marche}
F.~Marche, \emph{{Derivation of a new two-dimensional viscous shallow water
  model with varying topography, bottom friction and capillary effects}},
  European Journal of Mechanic /B \textbf{26} (2007), 49--63.

\bibitem{valentin}
B.~Mohammadi, O.~Pironneau, and F.~Valentin, \emph{{Rough boundaries and wall
  laws}}, Internat. J. Numer. Methods Fluids \textbf{27} (1998), no.~1-4,
  169--177.

\bibitem{nwogu}
O.~Nwogu, \emph{{Alternative form of Boussinesq equations for nearshore wave
  propagation}}, Journal of Waterway, Port, Coastal and Ocean Engineering, ASCE
  \textbf{119} (1993), no.~6, 618--638.

\bibitem{peregrine}
D.H. Peregrine, \emph{{Long waves on a beach}}, J. Fluid Mech. \textbf{27}
  (1967), 815--827.

\bibitem{perthame}
B.~Perthame, \emph{Kinetic formulation of conservation laws.}, Oxford
  University Press, 2002.

\bibitem{simeoni}
B.~Perthame and C.~Simeoni, \emph{A kinetic scheme for the saint-venant system
  with a source term}, Calcolo \textbf{38} (2001), no.~4, 201--231.

\bibitem{salencon}
M.J. Salen\c{c}on and J.M. Th\'ebault, \emph{Simulation model of a mesotrophic
  reservoir (lac de pareloup, france): Melodia, an ecosystem reservoir
  management model}, Ecological modelling \textbf{84} (1996), 163--187.

\bibitem{ursell}
F.~Ursell, \emph{{The long wave paradox in the theory of gavity waves}}, Proc.
  Cambridge Phil. Soc. \textbf{49} (1953), 685--694.

\bibitem{walkley}
M.A. Walkley, \emph{{A numerical Method for Extended Boussinesq Shallow-Water
  Wave Equations}}, Ph.D. thesis, University of Leeds, 1999.

\end{thebibliography}
